\documentclass[twoside,11pt,letterpaper]{article}
\usepackage[textwidth=160mm,textheight=220mm,includefoot,footskip=14mm,centering]{geometry}
\usepackage[english]{babel}
\usepackage[T1]{fontenc}
\usepackage{lmodern,microtype}
\usepackage{amssymb,amsmath,amsthm,indentfirst}
\usepackage[resetstyle,sepcounters,countpersection,proofenv=false,thmheading=false]{phfthm}
\usepackage[numbers,sort&compress]{natbib}
\usepackage{xcolor}
\usepackage[colorlinks=true,hyperfootnotes=false]{hyperref}
\hypersetup{
  linkcolor=red!80!black,
  citecolor=blue!90!black,
  urlcolor=blue!80!black,
  pdftitle={Global boundedness in a chemo-repulsion system with nonmonotone sensitivities of double-exponential growth},
  pdfauthor={Songzhi Li, Changchun Liu and Ming Mei}
}
\usepackage{doi}

\renewcommand{\labelenumi}{\textup{(\arabic{enumi})}}
\numberwithin{equation}{section}
\allowdisplaybreaks

\begin{document}
\renewcommand{\proofname}{\normalfont\bfseries\itshape Proof}
\title{\fontsize{16}{20}\selectfont\bfseries Global boundedness in a chemo-repulsion system with\\
nonmonotone sensitivities of double-exponential growth}
\author{
  Songzhi Li$^{1}$, Changchun Liu$^{2}$ and Ming Mei$^{1,3,4}$\\[1mm]
  \small\itshape $^1$School of Mathematics and Statistics, Jiangxi Normal University,\\
  \small\itshape Nanchang 330022, China\\
  \small\itshape $^2$School of Mathematics, Jilin University, Changchun 130012, China\\
  \small\itshape $^3$Department of Mathematics, Champlain College Saint-Lambert,\\
  \small\itshape Saint-Lambert, Quebec J4P 3P2, Canada\\
  \small\itshape $^4$Department of Mathematics and Statistics, McGill University,\\
  \small\itshape Montreal, Quebec H3A 2K6, Canada
}
\date{}
\maketitle
{\renewcommand{\thefootnote}{}\footnotetext{\raggedright E-mail addresses. lisz@jxnu.edu.cn (S. Li), liucc@jlu.edu.cn (C. Liu), ming.mei@mcgill.ca (M. Mei).}}

\vspace{-1.8\baselineskip}
{\setlength{\leftmargini}{1em}
\begin{quote}
\small
\textbf{Abstract.} This paper is concerned with a two-dimensional fully parabolic
chemo-repulsion system with linear diffusion, logistic growth and
signal-dependent sensitivity $\chi(v)$. The main issue is to prove the global existence of classic solution 
with uniform boundedness of the system for the challenging case with arbitrary admissible initial data when the sensitivity grows faster than any polynomial. However, the lack of pointwise control makes it difficult to estimate the sensitivity $\chi(v)$ for this challenging case. To overcome it by combining some technical energy estimates for exponential and double-exponential weights with entropy estimates, the key technical innovation developed in this paper is to
introduce a time-dependent shift $h(t)$ in the weight $w(v+h(t))$, with $h'(t)$ chosen
in terms of the entropy dissipation. The results significantly improve the existing studies, and   can be reviewed as the first framework of global bounedeness for arbitrary initial data with nonmonotone sensitivity of double-exponential growth. 

\indent\textbf{Keywords.} Chemo-repulsion, global classical solution, boundedness,
double-exponential sensitivity, logistic damping

\indent\textbf{AMS (2020) Subject Classification.} 35K51, 35A01, 35B45, 92C17
\end{quote}}

{\renewcommand{\addvspace}[1]{\vspace{1pt}}
\tableofcontents
\par}
\thispagestyle{empty}

\section{Introduction and main result}
\subsection{Model and motivation}\label{subsec-motivation}
In this paper, we study the following two-dimensional chemotaxis system
with linear diffusion, signal-dependent sensitivity and the logistic source $au-bu^2$
\begin{equation}\label{eq-model}
	\begin{cases}
		u_t=\Delta u+\nabla\cdot(u\chi(v)\nabla v)+au-bu^2,
		&x\in\Omega, t>0,\\
		v_t=\Delta v-v+u,&x\in\Omega, t>0,\\
		\partial_\nu u=\partial_\nu v=0,&x\in\partial\Omega, t>0,\\
		u(x,0)=u_0(x),\quad v(x,0)=v_0(x),&x\in\Omega.
	\end{cases}
\end{equation}
Here $u$ denotes the cell density and $v$ the concentration of the
chemical signal produced by the cells. With the sign convention in
\eqref{eq-model}, cells move toward higher signal concentrations when
$\chi(v)<0$ (chemo-attraction) and toward lower concentrations when
$\chi(v)>0$ (chemo-repulsion). In either case, $|\chi(v)|$ determines
the strength of the response to a given signal gradient. The domain
$\Omega\subset\mathbb R^2$ is bounded with smooth boundary, and $a,b>0$.

For constant sensitivity of either sign, the result of Osaki,
Tsujikawa, Yagi and Mimura \cite{OTYM} implies that in two dimensions
the logistic source $au-bu^2$ ensures global existence and uniform
boundedness for every $b>0$ without a smallness assumption on the initial data.
In the attractive case
$\chi\equiv-\chi_0<0$, Jin and Xiang \cite{JX2018} obtained bounds for
$\|u(\cdot,t)\|_{L^\infty(\Omega)}$ and
$\|v(\cdot,t)\|_{W^{1,\infty}(\Omega)}$ with explicit dependence on
$\chi_0$ and $b$.

The basic energy estimate for the signal equation bounds $v$ uniformly
in $H^1(\Omega)$. In two dimensions, $H^1(\Omega)$ does not embed into
$L^\infty(\Omega)$, so this estimate gives no uniform pointwise upper
bound for $v$. This lack of pointwise control makes it difficult to
estimate the factor $\chi(v)$ in the chemotactic term when the
sensitivity is allowed to grow without bound. To our knowledge, few results
establish global boundedness for \eqref{eq-model} with such sensitivities
and arbitrary initial data.

Recently, Li and Liu \cite{LL} treated sensitivities of at most polynomial
growth in the repulsive case by introducing the weighted energy
$\int_{\Omega}uv^{2p}+\int_{\Omega}|\nabla v^{p+1}|^2$ for $p\ge1$.
Differentiating a power of $v$ lowers its degree, so estimates for lower
powers can be used at the next step of the iteration. In this way they
controlled the chemotactic term and proved global boundedness for any
polynomial degree. It is then natural to ask whether this iteration can be
extended to sensitivities that grow faster than any polynomial. Let us consider first exponential growth, and 
differentiate the weight $e^{\lambda v}$ gives $\lambda e^{\lambda v}$,
then we recognize that  there is no reduction in the exponential rate to exploit at the next
step. The difficulty is greater still for double-exponential growth,
since differentiating a weight such as $e^{e^v}$ introduces the
additional unbounded factor $e^v$. Thus, this leads to the following important question:
{\it With the logistic damping fixed, can we still guarantee global,
uniformly bounded classical solutions for arbitrary admissible initial
data when the sensitivity is allowed to grow faster than any polynomial?}

In this paper, we address this question for chemo-repulsion, allowing the
sensitivity to be unbounded and nonmonotone.

\subsection{Research background}

\emph{Models with linear diffusion.}
The influence of the chemotactic term on boundedness can already be seen
in the classical Keller--Segel system \cite{KS,KSmodel}, obtained from
\eqref{eq-model} by taking $\chi$ constant and $a=b=0$. In the attractive case, Osaki and
Yagi \cite{OY} proved global existence and uniform boundedness in one
dimension for arbitrary initial mass. In two dimensions,
Nagai, Senba and Yoshida \cite[Theorem 1.1]{NSY} proved boundedness for
initial mass below $\frac{4\pi}{\chi_0}$, or below $\frac{8\pi}{\chi_0}$ for
radial data on a disk. Horstmann and Wang \cite{HW} constructed supercritical initial
data leading to blow-up in finite or infinite time. In three dimensions,
Winkler \cite{Wblowup} proved finite-time blow-up on a ball for suitable
radial initial data of arbitrarily small mass.

In the repulsive case $\chi\equiv\chi_0>0$, the mixed gradient terms cancel
upon differentiation of
$\int_{\Omega}u\log u+\frac{\chi_0}{2}\int_{\Omega}|\nabla v|^2$.
Cie\'slak, Lauren\c{c}ot and Morales-Rodrigo \cite{CLM} used this
cancellation to prove global boundedness and convergence to the homogeneous
equilibrium for arbitrary initial mass in two dimensions. They also
constructed global weak solutions
in dimensions three and four. In three dimensions, Jiang \cite{Jeventual}
proved the existence of global weak solutions that become smooth after
some finite time, and Cie\'slak, Fuest, Hajduk and Sier\.z\k{e}ga
\cite{CFHS} obtained a criterion for global classical solvability.
Whether the fully parabolic repulsive system without a source term admits
global classical solutions for arbitrary smooth initial data in dimensions
three and higher remains open.

For constant sensitivity, Winkler \cite[Theorem 0.1]{Wlogistic2010}
obtained a higher-dimensional boundedness result under sufficiently strong
logistic damping. In two dimensions, Winkler
\cite{WOrlicz} proved boundedness for every positive quadratic damping
coefficient when the sensitivity and its first derivative are bounded. With
his sign convention, the case $D\equiv1$ and $S=-\chi$ gives this result
for \eqref{eq-model}.

A classical example of signal dependence is logarithmic sensitivity, for
which cells respond to the relative signal gradient
$\nabla\log v=\frac{\nabla v}{v}$. For chemo-attraction with
$|\chi(v)|=\frac{\chi_0}{v}$, a strictly positive initial signal and no
source term, Winkler \cite{Wsingular}
proved global existence of classical solutions in dimensions $n\ge2$ for
$0<\chi_0<\sqrt{\frac{2}{n}}$, and Fujie \cite{Fsingular} proved uniform
boundedness in the same range. Mizukami and Yokota \cite{MYsingular}
obtained boundedness for stronger singularities of the form
$\frac{K}{v^k}$, $k>1$, under a smallness condition on $K$ depending on
the initial data. For chemo-repulsion with logarithmic sensitivity and
no source term, Jiang \cite{Jlog} used the energy
$E_0(\rho,c)=\int_{\Omega}(\rho\log\rho+2|\nabla\sqrt{c}|^2)$ for the
cell density $\rho$ and signal concentration $c$ to prove global classical
solvability in two dimensions for a strictly positive initial signal.
Chen, Kong and Wang \cite{CKWlog} proved exponential convergence to the
constant equilibrium in two dimensions for arbitrary repulsion strength.

For chemo-attraction with logarithmic sensitivity and logistic growth,
Zhao and Zheng \cite{ZZsingular} proved global boundedness in two
dimensions when the linear growth rate is sufficiently large relative
to the sensitivity coefficient. For the weaker singularity
$|\chi(v)|=\frac{\chi_0}{v^k}$, $0<k<1$, Le \cite{LeWeak} obtained
uniform boundedness of the cell density in two dimensions for every
positive quadratic damping coefficient. These inverse powers are
singular at zero and decay at high signal concentrations. In
\eqref{eq-model}, the sensitivity is regular at zero and may grow without
bound as $v\to\infty$.

\emph{Models with signal-dependent diffusion.}
Unbounded growth of the chemotactic sensitivity also arises in a widely
studied class of models with signal-dependent diffusion. In these models
the cell equation contains the term
\[
\Delta(u\gamma(v))=\nabla\cdot(\gamma(v)\nabla u)
+\nabla\cdot(u\gamma'(v)\nabla v).
\]
Thus the motility $\gamma(v)$ determines both the diffusion coefficient
and the chemotactic sensitivity $\gamma'(v)$. For example, the choice
$\gamma(v)=v^\alpha$ gives a repulsive sensitivity $\alpha v^{\alpha-1}$
that grows without bound when $\alpha>1$. For this family of motilities
and the generalized logistic source $au-bu^\sigma$, Li and Liu
\cite{LLmotility} proved global existence of classical solutions in
dimensions one and two for all $\alpha>0$ and $\sigma>1$. In the absence
of a source term, Jiang and Lauren\c{c}ot \cite{JL2026} proved that for
positive nondecreasing motilities, classical solutions of the fully
parabolic system are uniformly bounded in every dimension and converge
to the spatially homogeneous equilibrium. Their result imposes no upper
growth restriction on $\gamma'$.

Decreasing motilities correspond to chemo-attraction. In two dimensions,
Fujie and Jiang \cite{FJcomparison} showed that for $\gamma(v)=e^{-v}$
and no source term, solutions are bounded for subcritical mass and
can exhibit infinite-time blow-up for certain supercritical initial
data. For a class of decreasing motilities including $e^{-v}$, Jin,
Kim and Wang \cite{JKW2018} proved global existence and uniform
boundedness in two dimensions for every $\mu>0$ when the cell equation
contains the logistic source $\mu u(1-u)$.

Boundedness has also been studied for motilities that need not be
monotone. An early two-dimensional result for uniformly positive
motilities $\phi\in C^3([0,\infty))\cap W^{1,\infty}((0,\infty))$ is
due to Tao and Winkler \cite[Theorem 1.1]{TWmotility}. Xiao and Jiang \cite{XJ}
proved uniform boundedness in every dimension for motilities bounded
above and away from zero, without a global bound on $\gamma'$. Jiang and
Lauren\c{c}ot obtained boundedness for an elliptic signal equation under
$\limsup_{s\to\infty}\gamma(s)=\infty$ \cite{JL2024} and for a
parabolic signal equation under
$\liminf_{s\to\infty}\gamma(s)>\frac{1}{\tau}$ \cite{JL2026}, where
$\tau>0$ is the coefficient of $\partial_t v$ in the signal equation.

The product $u\gamma(v)$ plays a central role in the comparison methods
used for these models. For the source-free equation
$u_t=\Delta(u\gamma(v))$, Fujie and Jiang \cite{FJcomparison} introduced
$w=(I-\Delta)^{-1}[u]$ and obtained the identity
\[
\partial_t w+u\gamma(v)=(I-\Delta)^{-1}[u\gamma(v)].
\]
They used elliptic and parabolic comparison to bound $v$ pointwise in
terms of $w$. Xiao and Jiang \cite{XJ} and Jiang and Lauren\c{c}ot
\cite{JL2026} developed this approach for nonmonotone motilities. In
\eqref{eq-model}, matching the diffusion term $\Delta u$ would force
$\gamma\equiv1$ and hence $\gamma'\equiv0$. Unfortunately, system \eqref{eq-model}
lacks this product structure when $\chi\not\equiv0$, so this elegant
comparison method no longer applies.

\subsection{Main result and proof strategy}
Our main result gives a partial answer to the question raised in
Section \ref{subsec-motivation}. For any fixed $b>0$, we prove global
boundedness in the repulsive case for a class of sensitivities that need
not be monotone and may grow as fast as a double exponential.
To state the result, we assume that the initial data and the sensitivity satisfy
\begin{equation}\label{eq-data}
	0\le u_0\in C^0(\overline\Omega),\qquad
	0\le v_0\in W^{1,\infty}(\Omega),\qquad
	0\le\chi\in C^2([0,\infty)).
\end{equation}

We now state our main result.

\begin{theorem}\label{thm-main}
	Let $\Omega\subset\mathbb R^2$ be a bounded domain with smooth
	boundary. Assume \eqref{eq-data} and fix $a,b>0$. Suppose that at least
	one of the following conditions holds.
	\renewcommand{\labelenumi}{\textup{(\roman{enumi})}}
	\begin{enumerate}
		\item There exists $C>0$ such that
		\begin{equation}\label{eq-value-assumption}
			0\le\chi(s)\le Ce^{bs},\qquad(s\ge0).
		\end{equation}
		\item There exists $L>0$ such that
		\begin{equation}\label{eq-main-assumption}
			|\chi'(s)|\le L e^{\frac{bs}{2}}(1+\chi(s)),\qquad(s\ge0).
		\end{equation}
	\end{enumerate}
	Then problem \eqref{eq-model} has a unique global nonnegative
	classical solution, and
	\begin{equation}\label{eq-bounded}
		\sup_{t\ge0}\left(
		\|u(t)\|_{L^\infty(\Omega)}
		+\|v(t)\|_{W^{1,\infty}(\Omega)}\right)<\infty.
	\end{equation}
\end{theorem}

\begin{remark}
	Condition \eqref{eq-value-assumption} includes all sensitivities of at most polynomial
	growth considered by Li and Liu \cite{LL}, since $(1+s)^k\le C_{k,b}e^{bs}$
	for every $k\ge0$ and $s\ge0$.
\end{remark}

\begin{remark}
	Condition \eqref{eq-main-assumption} gives a double-exponential upper bound for $\chi$.
	Indeed, dividing \eqref{eq-main-assumption} by $1+\chi(s)$ and integrating from $0$ to $s$ gives
	\begin{equation}\label{eq-automatic-growth}
		1+\chi(s)\le(1+\chi(0))
		\exp\left\{\frac{2L}{b}\bigl(e^{\frac{bs}{2}}-1\bigr)\right\}.
	\end{equation}
\end{remark}

\begin{remark}\label{rem-examples}
	Theorem \ref{thm-main} applies to the following sensitivities.
	\begin{enumerate}
		\item For arbitrary $A,\beta,m>0$, the sensitivity $\chi(s)=A\exp\{\beta(1+s)^m\}$
		satisfies \eqref{eq-main-assumption}.
		\item Let $A,\beta>0$, $0<\gamma\le\frac{b}{2}$ and $x=e^{\gamma s}$.
		Both $\chi(s)=Ae^{\beta x}$ and $\chi(s)=Ae^{\beta(x+2\sin x)}$ satisfy
		\eqref{eq-main-assumption}, since $|\chi'(s)|\le3\beta\gamma e^{\gamma s}\chi(s)$.
		Both grow faster than $e^{bs}$, and the derivative of the second changes sign infinitely often.
		\item The oscillatory sensitivity $\chi(s)=e^{bs}[1+\sin(e^{bs})]$ satisfies
		\eqref{eq-value-assumption} with $C=2$. To check \eqref{eq-main-assumption}, take
		$s_n=\frac{\log(2\pi n)}{b}$ for $n\ge1$. Then $\chi(s_n)=2\pi n$ and
		$\chi'(s_n)=b(2\pi n)(1+2\pi n)$. Hence \eqref{eq-main-assumption} would require
		$b\sqrt{2\pi n}\le L$ for every $n\ge1$, which is impossible.
	\end{enumerate}
	The last two examples show that neither hypothesis of Theorem
	\ref{thm-main} implies the other.
\end{remark}

\begin{remark}
	The rates $b$ and $\frac{b}{2}$ enter at different points in the entropy
	estimates. Under \eqref{eq-value-assumption}, we use
	$\chi(v)^2\le Ce^{bv}\chi(v)$ and the bound for
	$\int_t^{t+\tau}\int_{\Omega}u\chi(v)e^{bv}|\nabla v|^2$
	from Lemma \ref{lem-low}. If $e^{\gamma s}$ replaces
	$e^{\frac{bs}{2}}$ in \eqref{eq-main-assumption}, the weighted
	interpolation step requires a bound for
	$\int_{\Omega}e^{2\gamma v}|\nabla v|^2$. Lemma \ref{lem-low}
	gives this bound for $2\gamma\le b$. Whether the value bound
	$\chi(s)\le Ce^{\kappa s}$ alone suffices for $\kappa>b$ is not
	settled by these estimates.
\end{remark}

\textbf{Difficulties and proof strategy.}
Following Li and Liu \cite{LL}, we seek a uniform $L^2(\Omega)$ bound
for $u$ by testing the cell equation with $u$. Two-dimensional
interpolation and the uniform Gr\"onwall inequality reduce this step to
the fourth-power estimate
\[
\int_t^{t+\tau}\int_{\Omega}\chi(v)^4|\nabla v|^4\le C,
\]
with $\tau$ given by \eqref{eq-tau} and $C$ independent of $t$.

We first consider the exponentially weighted energy
$\int_{\Omega}ue^{qv}+q\int_{\Omega}e^{qv}|\nabla v|^2$.
Signal production and logistic damping contribute the term
$(q-b)\int_{\Omega}u^2e^{qv}$ to its time derivative. For $0<q\le b$
this term is nonpositive. For $q>b$ it becomes positive, and the
estimates for smaller exponential weights do not control it.
To estimate this term, we apply the two-dimensional Sobolev inequality
to $ue^{\frac{qv}{2}}$. The resulting estimate contains
$\int_{\Omega}\frac{|\nabla u|^2}{u}$ as a coefficient of
$\int_{\Omega}ue^{qv}$, which is already part of the weighted energy.
The entropy estimate, obtained using the weights with $0<q\le b$,
bounds the integral of this coefficient over $(t,t+\tau)$ uniformly in $t$.
The uniform Gr\"onwall inequality then bounds the exponentially weighted energy for
every fixed $q>b$. Under \eqref{eq-value-assumption}, taking $q=2b$
yields the required fourth-power estimate.

Under \eqref{eq-main-assumption}, no fixed exponential weight covers
the double-exponential growth allowed by \eqref{eq-automatic-growth}.
We therefore take $w(s)=\exp\{Ae^{\frac{bs}{2}}\}$ with $A$ large
enough that $\chi^4\le Cw$, and consider the corresponding weighted energy
\[
\int_{\Omega}uw(v)+\int_{\Omega}w'(v)|\nabla v|^2.
\]
The mixed terms involving $\nabla u\cdot\nabla v$ cancel when we
differentiate this weighted energy.
Estimating the signal production term $\int_{\Omega}u^2w'(v)$ leaves
$\int_{\Omega}uw'(v)$ on the right-hand side of the weighted energy inequality.
Since $\frac{w'(s)}{w(s)}=\frac{Ab}{2}e^{\frac{bs}{2}}$ is unbounded,
the first term of the weighted energy gives no direct bound for
$\int_{\Omega}uw'(v)$.

The factor $w'(v)$ suggests shifting the argument of $w$.
Indeed, differentiating $\int_{\Omega}uw(v+h)$ in time produces
$h'\int_{\Omega}uw'(v+h)$, whose coefficient we can choose.
A shift changes only the parameter $A$ in the weight, while the
coefficient $C(1+\int_{\Omega}\frac{|\nabla u|^2}{u})$ of
$\int_{\Omega}uw'(v+h)$ in the estimate is independent of $A$.
Taking $h'$ to be the negative of this coefficient cancels the term.
With $h(t+\tau)=0$, the entropy dissipation bound keeps $h$ uniformly
bounded on $(t,t+\tau)$. We can thus close the weighted energy estimate
and obtain the required fourth-power bound by interpolation.

This paper is organized as follows. In Section \ref{sec-prelim}, we
recall the local existence result and present several estimates used in
the subsequent proofs. In Sections \ref{sec-exponential} and
\ref{sec-double-exponential}, we establish energy estimates with
exponential and double-exponential weights, respectively. Finally, in
Section \ref{sec-boundedness}, we use these estimates to derive a
space-time $L^4$ bound for $\chi(v)\nabla v$ and complete the proof of
Theorem \ref{thm-main}.

\section{Preliminaries}\label{sec-prelim}
We write $\int_{\Omega}f$ for $\int_{\Omega}fdx$ and
$\int_{t_1}^{t_2}\int_{\Omega}f$ for $\int_{t_1}^{t_2}\int_{\Omega}fdxds$.
The constants $K_i$ in the estimates and $c_i$ in the proofs are independent
of $t$. The numbering of the $c_i$ starts anew in each proof.

We begin by stating the local existence result for \eqref{eq-model} and
the corresponding extensibility criterion. See \cite{Amann,Amann1993} for the
parabolic theory and \cite[Lemma 2.2]{Cao2015} for a related chemotaxis system.

\begin{lemma}\label{lem-local}
Let $\Omega\subset\mathbb R^2$ be a bounded domain with smooth boundary,
let $a,b>0$ and assume \eqref{eq-data}. Then there exist a maximal
existence time $T_{\max}\in(0,\infty]$ and a unique nonnegative
classical solution $(u,v)$ of \eqref{eq-model} satisfying
\[
(u,v)\in\left[C^0(\overline\Omega\times[0,T_{\max}))
\cap C^{2,1}(\overline\Omega\times(0,T_{\max}))\right]^2.
\]
Moreover,
\begin{equation}\label{eq-alternative}
\text{if }T_{\max}<\infty,\quad\text{then}\quad
\limsup_{t\nearrow T_{\max}}
\left(\|u(\cdot,t)\|_{L^\infty(\Omega)}
+\|v(\cdot,t)\|_{W^{1,\infty}(\Omega)}\right)=\infty.
\end{equation}
\end{lemma}

If $u_0\equiv0$, then $u\equiv0$ and $v(\cdot,t)=e^{t(\Delta-1)}v_0$,
and Theorem \ref{thm-main} follows from the heat semigroup estimates
\eqref{eq-heat-Lp} and \eqref{eq-heat-initial-gradient}. In the sequel we assume
$u_0\not\equiv0$. The strong maximum principle and the Neumann boundary
condition then give $u>0$ on $\overline\Omega\times(0,T_{\max})$.

We state two forms of
Gr\"onwall's inequality. The first follows from the convolution estimate in
\cite[Lemma 3.4]{Wlogistic2019}.
\begin{lemma}\label{lem-damped-gronwall}
Let $T\in(0,\infty]$, $\tau\in(0,T)$ and $a,b>0$. Suppose that
$h\in L^1_{\rm loc}([0,T))$ is nonnegative and satisfies
\[
\int_t^{t+\tau}h(s)ds\le b,\quad\text{for all }t\in(0,T-\tau).
\]
Then
\[
\int_0^{t}e^{-a(t-s)}h(s)ds\le\frac{b}{1-e^{-a\tau}},
\quad\text{for all }t\in[0,T).
\]
Consequently, if $y\in C([0,T))$ is locally absolutely continuous on
$(0,T)$, $y(0)\ge0$ and $y'+ay\le h$ almost everywhere in $(0,T)$, then
\[
y(t)\le e^{-at}y(0)+\frac{b}{1-e^{-a\tau}}
\le y(0)+\frac{b}{1-e^{-a\tau}},
\quad\text{for all }t\in[0,T).
\]
\end{lemma}

The second is the uniform Gr\"onwall lemma
(see \cite[Lemma III.1.1]{Temam}).
\begin{lemma}\label{lem-uniform-gronwall}
Let $T\in(0,\infty]$ and $\tau\in(0,T)$.
Let $y$ be nonnegative and locally absolutely continuous on $(0,T)$,
and let $h,f\in L^1_{\rm loc}(0,T)$ be nonnegative. Suppose that
\[
y'(t)\le h(t)y(t)+f(t),
\quad\text{for a.e. }t\in(0,T).
\]
Assume that there are constants $A,B,D\ge0$ such that
\[
\int_t^{t+\tau}h(s)ds\le A,\qquad
\int_t^{t+\tau}f(s)ds\le B,\qquad
\int_t^{t+\tau}y(s)ds\le D,
\]
for all $t\in(0,T-\tau)$. Then
\[
y(t+\tau)\le \left(\frac{D}{\tau}+B\right)e^A,
\quad\text{for all }t\in(0,T-\tau).
\]
\end{lemma}
\begin{proof}
Fix $t\in(0,T-\tau)$. By the continuity of $y$ and the bound on its integral, there exists
$t_\star\in(t,t+\tau)$ such that
\[
y(t_\star)\le\frac{1}{\tau}\int_t^{t+\tau}y(s)ds\le\frac{D}{\tau}.
\]
Applying Gr\"onwall's inequality on $(t_\star,t+\tau)$, we obtain
\[
y(t+\tau)\le
\left(y(t_\star)+\int_{t_\star}^{t+\tau}f(s)ds\right)
\exp\left(\int_{t_\star}^{t+\tau}h(s)ds\right)
\le \left(\frac{D}{\tau}+B\right)e^A.
\]
This completes the proof.
\end{proof}

\filbreak
We now derive the basic estimates for $u$ and $v$.
\begin{lemma}\label{lem-basic}
Let $(u,v)$ be the solution given by Lemma \ref{lem-local}. There exist
$K_i>0$ $(i=1,2,3,4)$, independent of $t$, such that
\begin{equation}\label{eq-mass}
\int_{\Omega}u(\cdot,t)\le K_1,
\quad\text{for all }t\in(0,T_{\max}),
\end{equation}
and
\begin{equation}\label{eq-u2-window}
\int_t^{t+\tau}\int_{\Omega}u^2\le K_2,
\quad\text{for all }t\in(0,\widetilde T_{\max}).
\end{equation}
Moreover,
\begin{equation}\label{eq-signal-H1}
\int_{\Omega}v^2(\cdot,t)+\int_{\Omega}|\nabla v(\cdot,t)|^2\le K_3,
\quad\text{for all }t\in(0,T_{\max}),
\end{equation}
and
\begin{equation}\label{eq-laplace-window}
\int_t^{t+\tau}\int_{\Omega}|\Delta v|^2\le K_4,
\quad\text{for all }t\in(0,\widetilde T_{\max}),
\end{equation}
where
\begin{equation}\label{eq-tau}
\tau=\min\left\{1,\frac{1}{2}T_{\max}\right\},\qquad
\widetilde T_{\max}=
\begin{cases}
T_{\max}-\tau,&\text{if }T_{\max}<\infty,\\
\infty,&\text{if }T_{\max}=\infty.
\end{cases}
\end{equation}
\end{lemma}
\begin{proof}
Integrating the first equation of \eqref{eq-model} over $\Omega$ and
using the boundary conditions, we obtain
\[
\frac{d}{dt}\int_{\Omega}u+b\int_{\Omega}u^2=a\int_{\Omega}u.
\]
By Young's inequality,
\[
\frac{d}{dt}\int_{\Omega}u+\int_{\Omega}u
\le\frac{(a+1)^2}{4b}|\Omega|,
\]
which proves \eqref{eq-mass} by Gr\"onwall's inequality. Integrating the first
equation of \eqref{eq-model} over $\Omega\times(t,t+\tau)$ and using \eqref{eq-mass}, we find
\[
b\int_t^{t+\tau}\int_{\Omega}u^2\le (1+a\tau)K_1,
\]
which proves \eqref{eq-u2-window}.

We test the second equation of \eqref{eq-model} with $v$ and $-\Delta v$
and add the two identities. Estimating $\int_{\Omega}uv$ and
$-\int_{\Omega}u\Delta v$ by Young's inequality, we obtain
\begin{equation}\label{eq-signal-differential}
\frac{d}{dt}\left(\int_{\Omega}v^2+\int_{\Omega}|\nabla v|^2\right)
+\int_{\Omega}v^2+\int_{\Omega}|\nabla v|^2+\int_{\Omega}|\Delta v|^2
\le 2\int_{\Omega}u^2.
\end{equation}
Using \eqref{eq-u2-window}, we apply Lemma \ref{lem-damped-gronwall} to
\eqref{eq-signal-differential} to obtain \eqref{eq-signal-H1}. Integrating
\eqref{eq-signal-differential} over $(t,t+\tau)$ and using \eqref{eq-u2-window}
and \eqref{eq-signal-H1}, we obtain \eqref{eq-laplace-window}.
This completes the proof.
\end{proof}

We shall also use the following two-dimensional Gagliardo--Nirenberg inequality
(see \cite[Lemma 2.5]{JKW2018}).
\begin{lemma}\label{lem-GN}
Let $\Omega$ be a bounded domain in $\mathbb R^2$ with smooth
boundary. There exists $C>0$, depending only on $\Omega$, such that
every $\varphi\in W^{2,2}(\Omega)$ with $\partial_\nu\varphi=0$ on
$\partial\Omega$ satisfies
\begin{equation}\label{eq-GN}
\|\nabla\varphi\|_{L^4(\Omega)}
\le C\left(\|\Delta\varphi\|_{L^2(\Omega)}^{\frac{1}{2}}
\|\nabla\varphi\|_{L^2(\Omega)}^{\frac{1}{2}}
+\|\nabla\varphi\|_{L^2(\Omega)}\right).
\end{equation}
\end{lemma}

The preceding estimates control $\int_t^{t+\tau}\int_{\Omega}|\nabla v|^4$
and give the first bounds involving exponential weights.
\begin{lemma}\label{lem-initial-exponential}
Let $(u,v)$ be the solution given by Lemma \ref{lem-local}.
There exists $K_5>0$, independent of $t$, such that
\begin{equation}\label{eq-gradient4-window}
\int_t^{t+\tau}\int_{\Omega}|\nabla v|^4\le K_5,
\quad\text{for all }t\in(0,\widetilde T_{\max}).
\end{equation}
For every $\lambda>0$, there also exists $K_6(\lambda)>0$ such that
\begin{equation}\label{eq-expmoment}
\int_{\Omega}e^{\lambda v(\cdot,t)}\le K_6(\lambda),
\quad\text{for all }t\in(0,T_{\max}),
\end{equation}
and
\begin{equation}\label{eq-expaverage}
\int_t^{t+\tau}\int_{\Omega}e^{\lambda v}|\nabla v|^2+\int_t^{t+\tau}\int_{\Omega}ue^{\lambda v}
\le K_6(\lambda),
\quad\text{for all }t\in(0,\widetilde T_{\max}).
\end{equation}
\end{lemma}
\begin{proof}
With $\varphi=v$ in \eqref{eq-GN}, we obtain
\[
\int_{\Omega}|\nabla v|^4
\le c_1\left(\int_{\Omega}|\nabla v|^2\right)
\left(\int_{\Omega}|\Delta v|^2\right)
+c_1\left(\int_{\Omega}|\nabla v|^2\right)^2.
\]
To prove \eqref{eq-gradient4-window}, we integrate over $(t,t+\tau)$
and use \eqref{eq-signal-H1} and \eqref{eq-laplace-window}.

By the Moser--Trudinger inequality \cite{Moser}, for each $M>0$ there are
constants $c_2,c_3>0$, depending only on $M$ and $\Omega$, such that
\begin{equation}\label{eq-TM}
\|z\|_{H^1(\Omega)}\le M
\Longrightarrow
\int_{\Omega}e^{c_2 z^2}\le c_3,
\end{equation}
for every $z\in H^1(\Omega)$.
In view of \eqref{eq-signal-H1}, we take $z=v$ in \eqref{eq-TM} and apply
Young's inequality to $\lambda v$ to obtain \eqref{eq-expmoment} for every $\lambda>0$.

To prove \eqref{eq-expaverage}, we first apply H\"older's inequality,
\eqref{eq-u2-window} and \eqref{eq-expmoment} with $2\lambda$ to obtain
\begin{equation}\label{eq-cell-exponential-average}
\int_t^{t+\tau}\int_{\Omega}ue^{\lambda v}
\le\left(\int_t^{t+\tau}\int_{\Omega}u^2\right)^{\frac{1}{2}}
\left(\int_t^{t+\tau}\int_{\Omega}e^{2\lambda v}\right)^{\frac{1}{2}}
\le c_4(\lambda).
\end{equation}
Next, multiplying the second equation of \eqref{eq-model} by $e^{\lambda v}$
and integrating by parts, we have
\begin{align*}
\frac{1}{\lambda}\frac{d}{dt}\int_{\Omega}e^{\lambda v}
=&\int_{\Omega}e^{\lambda v}(\Delta v-v+u)\\
=&-\lambda\int_{\Omega}e^{\lambda v}|\nabla v|^2
-\int_{\Omega}ve^{\lambda v}+\int_{\Omega}ue^{\lambda v}.
\end{align*}
Since $v\ge0$, integration over $(t,t+\tau)$ gives
\[
\lambda\int_t^{t+\tau}\int_{\Omega}e^{\lambda v}|\nabla v|^2
\le\frac{1}{\lambda}\int_{\Omega}e^{\lambda v(t)}
+\int_t^{t+\tau}\int_{\Omega}ue^{\lambda v}
\le c_5(\lambda).
\]
The last inequality follows from \eqref{eq-expmoment} and \eqref{eq-cell-exponential-average},
which also bounds the other integral in \eqref{eq-expaverage}.
\end{proof}

The next interpolation inequality estimates a weighted integral of $|\nabla v|^4$.
\begin{lemma}\label{lem-variable}
Let $(u,v)$ be the local classical solution given by Lemma \ref{lem-local}
and fix $t\in(0,T_{\max})$.
Let $\rho,R\in C^1([0,\infty))$ be positive, and suppose that there
exist $c,K>0$ such that
\[
\frac{|\rho'(s)|}{\rho(s)}\le cR(s),\qquad
\frac{|R'(s)|}{R(s)}\le cR(s)\quad(s\ge0),\qquad
\int_{\Omega}R(v)^2|\nabla v|^2\le K.
\]
Then for every $\varepsilon>0$ there exists
$K_7(\varepsilon)>0$, depending only on $\varepsilon,\Omega,c$ and $K$, such that
\begin{align}\label{eq-variable}
\int_{\Omega}\rho(v)R(v)^2|\nabla v|^4
\le&\varepsilon\int_{\Omega}\rho(v)|\Delta v|^2\notag\\
&+K_7(\varepsilon)\left(1+\int_{\Omega}R(v)^2|\Delta v|^2
+\int_{\Omega}R(v)^4|\nabla v|^4\right)
\int_{\Omega}\rho(v)|\nabla v|^2.
\end{align}
\end{lemma}
\begin{proof}
	To apply \eqref{eq-GN}, set
	\[
	\psi(s)=\int_0^s\rho(\sigma)^{\frac{1}{4}}R(\sigma)^{\frac{1}{2}}d\sigma,
	\qquad s\ge0,
	\]
	so that $(\psi')^4=\rho R^2$. Differentiating $\psi'$ and using the
	hypotheses on $\rho'$ and $R'$, we obtain
	\begin{align*}
		\psi''(s)
		=&\frac{1}{4}\rho(s)^{-\frac{3}{4}}\rho'(s)R(s)^{\frac{1}{2}}
		  +\frac{1}{2}\rho(s)^{\frac{1}{4}}R(s)^{-\frac{1}{2}}R'(s),\\
		|\psi''(s)|
		\le& \rho(s)^{\frac{1}{4}}R(s)^{\frac{1}{2}}
		     \left(\frac{|\rho'(s)|}{4\rho(s)}
		               +\frac{|R'(s)|}{2R(s)}\right)\\
		\le& \rho(s)^{\frac{1}{4}}R(s)^{\frac{1}{2}}
		     \left(\frac{c}{4} R(s)+\frac{c}{2} R(s)\right)
		 =\frac{3c}{4}\rho(s)^{\frac{1}{4}}R(s)^{\frac{3}{2}}.
	\end{align*}
	Applying \eqref{eq-GN} to $\psi(v)$, we find
	\begin{align*}
		\int_{\Omega}|\nabla\psi(v)|^4
		\le& c_1^4\left(
		    \|\Delta\psi(v)\|_{L^2(\Omega)}^{\frac{1}{2}}
		    \|\nabla\psi(v)\|_{L^2(\Omega)}^{\frac{1}{2}}
		    +\|\nabla\psi(v)\|_{L^2(\Omega)}\right)^4\\
		\le& 8c_1^4\left(\int_{\Omega}|\nabla\psi(v)|^2\right)
		    \left(\int_{\Omega}|\Delta\psi(v)|^2
		          +\int_{\Omega}|\nabla\psi(v)|^2\right).
	\end{align*}
	Substituting $\nabla\psi(v)=\psi'(v)\nabla v$ and
	$\Delta\psi(v)=\psi'(v)\Delta v+\psi''(v)|\nabla v|^2$, we find
	\begin{align*}
		\int_{\Omega}\rho(v)R(v)^2|\nabla v|^4
		=&\int_{\Omega}|\nabla\psi(v)|^4\\
		\le& c_2\left(\int_{\Omega}\sqrt{\rho(v)}R(v)|\nabla v|^2\right)
		\left(\int_{\Omega}\sqrt{\rho(v)}R(v)|\Delta v|^2\right)\\
		&+c_2\left(\int_{\Omega}\sqrt{\rho(v)}R(v)|\nabla v|^2\right)
		\left(\int_{\Omega}\sqrt{\rho(v)}R(v)^3|\nabla v|^4\right)\\
		&+c_2\left(\int_{\Omega}\sqrt{\rho(v)}R(v)|\nabla v|^2\right)^2.
	\end{align*}
	By H\"older's inequality and the assumed bound on
	$\int_{\Omega}R(v)^2|\nabla v|^2$, the factors on the right satisfy
	\begin{align*}
		\int_{\Omega}\sqrt{\rho(v)}R(v)|\nabla v|^2
		\le& K^{\frac{1}{2}}\left(\int_{\Omega}\rho(v)|\nabla v|^2\right)^{\frac{1}{2}},\\
		\int_{\Omega}\sqrt{\rho(v)}R(v)|\Delta v|^2
		\le&\left(\int_{\Omega}\rho(v)|\Delta v|^2\right)^{\frac{1}{2}}
		\left(\int_{\Omega}R(v)^2|\Delta v|^2\right)^{\frac{1}{2}},\\
		\int_{\Omega}\sqrt{\rho(v)}R(v)^3|\nabla v|^4
		\le&\left(\int_{\Omega}\rho(v)R(v)^2|\nabla v|^4\right)^{\frac{1}{2}}
		\left(\int_{\Omega}R(v)^4|\nabla v|^4\right)^{\frac{1}{2}}.
	\end{align*}
	Substituting these estimates and applying Young's inequality to the
	two products containing $\Delta v$ and $|\nabla v|^4$, we obtain
	\begin{align*}
		\int_{\Omega}\rho(v)R(v)^2|\nabla v|^4
		\le&\frac{\varepsilon}{2}\int_{\Omega}\rho(v)|\Delta v|^2
		+\frac{1}{2}\int_{\Omega}\rho(v)R(v)^2|\nabla v|^4\\
		&+c_3(\varepsilon)\left(1+\int_{\Omega}R(v)^2|\Delta v|^2
		+\int_{\Omega}R(v)^4|\nabla v|^4\right)
		\int_{\Omega}\rho(v)|\nabla v|^2.
	\end{align*}
	This proves \eqref{eq-variable} with $K_7(\varepsilon)=2c_3(\varepsilon)$.
\end{proof}

The final lemma collects the Neumann heat semigroup estimates needed in the proof
of Theorem \ref{thm-main} (see \cite[Lemma 1.3]{Wsemigroup} and \cite[Lemma 2.1]{Cao2015}).

\begin{lemma}
Let $(e^{t\Delta})_{t\ge0}$ be the Neumann heat semigroup on $\Omega$,
and let $\lambda_1>0$ be the first nonzero eigenvalue of $-\Delta$
in $\Omega$ under Neumann boundary conditions. There exist constants
$\gamma_1,\ldots,\gamma_4>0$, depending only on $\Omega$, such that
the following estimates hold for all $t>0$.
\renewcommand{\labelenumi}{\textup{(\roman{enumi})}}
\begin{enumerate}
\item If $1\le q\le p\le\infty$, then
\begin{equation}\label{eq-heat-Lp}
\|e^{t\Delta}w\|_{L^p(\Omega)}
\le\gamma_1\left(1+t^{-\left(\frac{1}{q}-\frac{1}{p}\right)}\right)
\|w\|_{L^q(\Omega)}
\end{equation}
for all $w\in L^q(\Omega)$.
\item If $1\le q\le p\le\infty$, then
\begin{equation}\label{eq-heat-gradient}
\|\nabla e^{t\Delta}w\|_{L^p(\Omega)}
\le\gamma_2\left(1+t^{-\frac{1}{2}-\left(\frac{1}{q}-\frac{1}{p}\right)}\right)
e^{-\lambda_1t}\|w\|_{L^q(\Omega)}
\end{equation}
for all $w\in L^q(\Omega)$.
\item If $2\le p\le\infty$, then
\begin{equation}\label{eq-heat-initial-gradient}
\|\nabla e^{t\Delta}w\|_{L^p(\Omega)}
\le\gamma_3 e^{-\lambda_1t}\|\nabla w\|_{L^p(\Omega)}
\end{equation}
for all $w\in W^{1,p}(\Omega)$.
\item If $1<q<\infty$ and $q\le p\le\infty$, then
\begin{equation}\label{eq-heat-divergence}
\|e^{t\Delta}\nabla\cdot F\|_{L^p(\Omega)}
\le\gamma_4\left(1+t^{-\frac{1}{2}-\left(\frac{1}{q}-\frac{1}{p}\right)}\right)
e^{-\lambda_1t}\|F\|_{L^q(\Omega)}
\end{equation}
for all $F\in C^1(\overline\Omega;\mathbb R^2)$ with $F\cdot\nu=0$ on $\partial\Omega$.
\end{enumerate}
\end{lemma}

\section{Estimates with exponential weights}\label{sec-exponential}
Motivated by the algebraically weighted energy in \cite{LL}, we aim to bound
\[
\int_{\Omega}ue^{\lambda v}
+\lambda\int_{\Omega}e^{\lambda v}|\nabla v|^2
\]
uniformly in time for each fixed $\lambda>0$. We begin with $0<\lambda\le b$.

\begin{lemma}\label{lem-low}
For every $0<\lambda\le b$, there exists $K_8(\lambda)>0$ such that
\begin{equation}\label{eq-low-pointwise}
\int_{\Omega}ue^{\lambda v}
+\int_{\Omega}e^{\lambda v}|\nabla v|^2\le K_8(\lambda),
\quad\text{for all }t\in(0,T_{\max}).
\end{equation}
Moreover,
\begin{equation}\label{eq-low-dissipation}
\int_t^{t+\tau}\int_{\Omega}e^{\lambda v}|\Delta v|^2
+\int_t^{t+\tau}\int_{\Omega}u\chi(v)e^{\lambda v}|\nabla v|^2
+\int_t^{t+\tau}\int_{\Omega}e^{2\lambda v}|\nabla v|^4\le K_8(\lambda),
\end{equation}
for all $t\in(0,\widetilde T_{\max})$.
\end{lemma}
\begin{proof}
Let $\lambda>0$. Differentiating $\int_{\Omega}ue^{\lambda v}$ in time,
substituting from \eqref{eq-model} and integrating by parts, we find
\begin{align}\label{3.3}
\frac{d}{dt}\int_{\Omega}ue^{\lambda v}
=&\int_{\Omega}e^{\lambda v}
[\Delta u+\nabla\cdot(u\chi(v)\nabla v)+au-bu^2]\nonumber\\
&+\lambda\int_{\Omega}ue^{\lambda v}(\Delta v-v+u)\nonumber\\
=&-\lambda\int_{\Omega}e^{\lambda v}\nabla u\cdot\nabla v
-\lambda\int_{\Omega}u\chi(v)e^{\lambda v}|\nabla v|^2
+\lambda\int_{\Omega}ue^{\lambda v}\Delta v\nonumber\\
&+(\lambda-b)\int_{\Omega}u^2e^{\lambda v}
+\int_{\Omega}ue^{\lambda v}(a-\lambda v)\nonumber\\
=&-2\lambda\int_{\Omega}e^{\lambda v}\nabla u\cdot\nabla v
-\lambda^2\int_{\Omega}ue^{\lambda v}|\nabla v|^2
-\lambda\int_{\Omega}u\chi(v)e^{\lambda v}|\nabla v|^2\nonumber\\
&+(\lambda-b)\int_{\Omega}u^2e^{\lambda v}
+\int_{\Omega}ue^{\lambda v}(a-\lambda v).
\end{align}
Using the second equation of \eqref{eq-model} and integration by
parts, we have
\begin{align}\label{3.4}
\frac{d}{dt}\int_{\Omega}e^{\lambda v}|\nabla v|^2
=&\lambda\int_{\Omega}e^{\lambda v}v_t|\nabla v|^2
+2\int_{\Omega}e^{\lambda v}\nabla v\cdot\nabla v_t\nonumber\\
=&-2\int_{\Omega}e^{\lambda v}\Delta vv_t
-\lambda\int_{\Omega}e^{\lambda v}|\nabla v|^2v_t\nonumber\\
=&-2\int_{\Omega}e^{\lambda v}\Delta v(\Delta v-v+u)
-\lambda\int_{\Omega}e^{\lambda v}|\nabla v|^2(\Delta v-v+u)\nonumber\\
=&-2\int_{\Omega}e^{\lambda v}|\Delta v|^2
-\lambda\int_{\Omega}e^{\lambda v}|\nabla v|^2\Delta v
+2\int_{\Omega}ve^{\lambda v}\Delta v\nonumber\\
&-2\int_{\Omega}ue^{\lambda v}\Delta v
+\lambda\int_{\Omega}ve^{\lambda v}|\nabla v|^2
-\lambda\int_{\Omega}ue^{\lambda v}|\nabla v|^2.
\end{align}
Integrating by parts once more in the two terms below, we find
\begin{align*}
2\int_{\Omega}ve^{\lambda v}\Delta v
=-2\int_{\Omega}e^{\lambda v}|\nabla v|^2-2\lambda\int_{\Omega}ve^{\lambda v}|\nabla v|^2,
\end{align*}
and
\begin{align*}
-2\int_{\Omega}ue^{\lambda v}\Delta v=2\int_{\Omega}e^{\lambda v}\nabla u\cdot\nabla v
 +2\lambda\int_{\Omega}ue^{\lambda v}|\nabla v|^2.
\end{align*}
Substituting the above two equations into \eqref{3.4} and multiplying by $\lambda$ yields
\begin{align}\label{3.5}
	&\lambda\frac{d}{dt}\int_{\Omega}e^{\lambda v}|\nabla v|^2+2\lambda\int_{\Omega}e^{\lambda v}|\Delta v|^2+2\lambda\int_{\Omega}e^{\lambda v}|\nabla v|^2+\lambda^2\int_{\Omega}ve^{\lambda v}|\nabla v|^2\nonumber\\*
	=&
	-\lambda^2\int_{\Omega}e^{\lambda v}|\nabla v|^2\Delta v
	+2\lambda\int_{\Omega}e^{\lambda v}\nabla u\cdot\nabla v
	+\lambda^2\int_{\Omega}ue^{\lambda v}|\nabla v|^2.
\end{align}
When we add \eqref{3.3} and \eqref{3.5}, the terms $2\lambda\int_{\Omega}e^{\lambda v}\nabla u\cdot\nabla v$ and $\lambda^2\int_{\Omega}ue^{\lambda v}|\nabla v|^2$ cancel, and we obtain
\begin{align}\label{eq-low-identity}
&\frac{d}{dt}\left(\int_{\Omega}ue^{\lambda v}
+\lambda\int_{\Omega}e^{\lambda v}|\nabla v|^2\right)+2\lambda\int_{\Omega}e^{\lambda v}|\Delta v|^2+\lambda^2\int_{\Omega}ve^{\lambda v}|\nabla v|^2+2\lambda\int_{\Omega}e^{\lambda v}|\nabla v|^2\notag\\
=&-\lambda^2\int_{\Omega}e^{\lambda v}|\nabla v|^2\Delta v-\lambda\int_{\Omega}u\chi(v)e^{\lambda v}|\nabla v|^2
+(\lambda-b)\int_{\Omega}u^2e^{\lambda v}+\int_{\Omega}ue^{\lambda v}(a-\lambda v).
\end{align}
By Young's inequality and Lemma \ref{lem-variable} with
$\rho(s)=e^{\lambda s}$ and $R=1$, using \eqref{eq-signal-H1} and
choosing $\varepsilon$ sufficiently small, we conclude that
\begin{align*}
	\lambda^2\left|\int_{\Omega}e^{\lambda v}|\nabla v|^2\Delta v\right|
\le&\frac{\lambda}{2}\int_{\Omega}e^{\lambda v}|\Delta v|^2
+c_1(\lambda)\int_{\Omega}e^{\lambda v}|\nabla v|^4\nonumber\\
\le& \lambda\int_{\Omega}e^{\lambda v}|\Delta v|^2+c_2(\lambda)\left(1+\int_{\Omega}|\Delta v|^2
+\int_{\Omega}|\nabla v|^4\right)
\int_{\Omega}e^{\lambda v}|\nabla v|^2.
\end{align*}
Since $s\mapsto e^{\lambda s}(a-\lambda s)$ is bounded above on
$[0,\infty)$, \eqref{eq-mass} gives
\begin{align*}
\int_{\Omega}ue^{\lambda v}(a-\lambda v)
\le c_3(\lambda)\int_{\Omega}u\le c_3(\lambda)K_1.
\end{align*}
Since $0<\lambda\le b$ and $\chi\ge0$, we have
\begin{align*}
	-\lambda\int_{\Omega}u\chi(v)e^{\lambda v}|\nabla v|^2+(\lambda-b)\int_{\Omega}u^2e^{\lambda v}\le-\lambda\int_{\Omega}u\chi(v)e^{\lambda v}|\nabla v|^2\le0.
\end{align*}
Substituting these estimates into \eqref{eq-low-identity} gives
\begin{align}\label{eq-lowineq}
&\frac{d}{dt}\left(\int_{\Omega}ue^{\lambda v}
+\lambda\int_{\Omega}e^{\lambda v}|\nabla v|^2\right)+\lambda\int_{\Omega}e^{\lambda v}|\Delta v|^2
+\lambda\int_{\Omega}u\chi(v)e^{\lambda v}|\nabla v|^2\notag\\
\le& \frac{c_2(\lambda)}{\lambda}\left(1+\int_{\Omega}|\Delta v|^2
+\int_{\Omega}|\nabla v|^4\right)
\left(\int_{\Omega}ue^{\lambda v}
+\lambda\int_{\Omega}e^{\lambda v}|\nabla v|^2\right)+c_3(\lambda)K_1.
\end{align}
By \eqref{eq-laplace-window}, \eqref{eq-gradient4-window}
and \eqref{eq-expaverage}, we have
\begin{align*}
\int_t^{t+\tau}\left(1+\int_{\Omega}|\Delta v|^2
+\int_{\Omega}|\nabla v|^4\right)
\le \tau+K_4+K_5,
\end{align*}
and
\begin{align*}
\int_t^{t+\tau}\int_{\Omega}ue^{\lambda v}+\int_t^{t+\tau}\int_{\Omega}\lambda e^{\lambda v}|\nabla v|^2
\le \max\{1,\lambda\}K_6(\lambda).
\end{align*}
Applying Lemma \ref{lem-uniform-gronwall} to \eqref{eq-lowineq} gives
\[
\int_{\Omega}ue^{\lambda v}
+\lambda\int_{\Omega}e^{\lambda v}|\nabla v|^2\le c_4(\lambda),
\qquad (\tau<t<T_{\max}).
\]
Since $v$ is bounded on $[0,\tau]$ by Lemma \ref{lem-local}, the preceding estimate
together with \eqref{eq-mass} and \eqref{eq-signal-H1} gives \eqref{eq-low-pointwise}
with upper bound $c_5(\lambda)$.

To bound the first two integrals in \eqref{eq-low-dissipation}, we integrate
\eqref{eq-lowineq} over $(t,t+\tau)$ and use \eqref{eq-low-pointwise},
\eqref{eq-laplace-window} and \eqref{eq-gradient4-window}.
For the last integral, we apply \eqref{eq-GN} to $e^{\frac{\lambda}{2}v}$
and then use \eqref{eq-low-pointwise} and Young's inequality to find
\begin{align}\label{eq-exponential-quartic}
	\int_{\Omega}e^{2\lambda v}|\nabla v|^4=& \frac{16}{\lambda^4}\int_{\Omega}\left|\nabla e^{\frac{\lambda}{2}v}\right|^4\nonumber\\
	\le& c_6(\lambda)\int_{\Omega}\left|\Delta e^{\frac{\lambda}{2}v}\right|^2\int_{\Omega}\left|\nabla e^{\frac{\lambda}{2}v}\right|^2
	+c_6(\lambda)\left(\int_{\Omega}\left|\nabla e^{\frac{\lambda}{2}v}\right|^2\right)^2\nonumber\\
	\le&c_7(\lambda)\left(\int_{\Omega}e^{\lambda v}|\nabla v|^2\right)
	\left(\int_{\Omega}e^{\lambda v}|\Delta v|^2
	+\int_{\Omega}e^{\lambda v}|\nabla v|^4\right)+c_7(\lambda)\left(\int_{\Omega}e^{\lambda v}|\nabla v|^2\right)^2\nonumber\\
	\le& c_8(\lambda)\left(1+\int_{\Omega}e^{\lambda v}|\Delta v|^2+\int_{\Omega}e^{\lambda v}|\nabla v|^4\right)\nonumber\\
	\le& \frac{1}{2}\int_{\Omega}e^{2\lambda v}|\nabla v|^4+ c_9(\lambda)\left(1+\int_{\Omega}e^{\lambda v}|\Delta v|^2+\int_{\Omega}|\nabla v|^4\right).
\end{align}
Integrating \eqref{eq-exponential-quartic} over $(t,t+\tau)$, we obtain
\begin{align*}
	\int_{t}^{t+\tau}\int_{\Omega}e^{2\lambda v}|\nabla v|^4\le 2c_9(\lambda)\int_{t}^{t+\tau}\left(1+\int_{\Omega}e^{\lambda v}|\Delta v|^2+\int_{\Omega}|\nabla v|^4\right).
\end{align*}
Thus, the bound for the last term in \eqref{eq-low-dissipation} follows from \eqref{eq-gradient4-window} and the bound for the first term in \eqref{eq-low-dissipation}. Choosing $K_8(\lambda)$ to cover both bounds completes the proof.
\end{proof}

Taking $\lambda=b$, we can now apply the Moser--Trudinger inequality
to $e^{\frac{bv}{2}}$. We shall use the resulting integral bounds in the
entropy estimate under \eqref{eq-main-assumption} and again in Section
\ref{sec-double-exponential}.
\begin{lemma}\label{lem-moments}
For every fixed $c>0$ and $0<\eta<b$, there exists
$K_9(c,\eta)>0$ such that
\begin{equation}\label{eq-base-double-moment}
\int_{\Omega}\exp\{ce^{\eta v}\}\le K_9(c,\eta),
\qquad\text{for all }t\in(0,T_{\max}).
\end{equation}
Moreover,
\begin{equation}\label{eq-weightedaverage}
\int_t^{t+\tau}\int_{\Omega}\exp\{ce^{\eta v}\}|\nabla v|^2\le K_9(c,\eta),
\qquad\text{for all }t\in(0,\widetilde T_{\max}).
\end{equation}
\end{lemma}
\begin{proof}
Taking $\lambda=b$ in \eqref{eq-expmoment} and \eqref{eq-low-pointwise}, we obtain
\[
\left\|e^{\frac{bv}{2}}\right\|_{H^1(\Omega)}^2
=\int_{\Omega}e^{bv}+\frac{b^2}{4}\int_{\Omega}e^{bv}|\nabla v|^2\le K_6(b)+\frac{b^2}{4}K_8(b).
\]
By the Moser--Trudinger inequality \eqref{eq-TM}, there exist $c_1(b),c_2(b)>0$ such that
\[
\int_{\Omega}\exp\{c_1(b)e^{bv}\}\le c_2(b),
\qquad\text{for all }t\in(0,T_{\max}).
\]
Since $0<\eta<b$, we apply Young's inequality with exponent $\frac{b}{\eta}$ to
$ce^{\eta v}$ to obtain, for each $c>0$ and some $c_3(c,\eta)>0$,
\[
\int_{\Omega}\exp\{ce^{\eta v}\}
\le e^{c_3(c,\eta)}\int_{\Omega}\exp\{c_1(b)e^{bv}\}\le c_4(c,\eta),
\]
which proves \eqref{eq-base-double-moment}.

For the second assertion, we apply H\"older's inequality in space and time to obtain
\begin{align*}
\int_t^{t+\tau}\int_{\Omega}\exp\{ce^{\eta v}\}|\nabla v|^2
\le&
\left(\int_t^{t+\tau}\int_{\Omega}\exp\{2ce^{\eta v}\}\right)^{\frac{1}{2}}
\left(\int_t^{t+\tau}\int_{\Omega}|\nabla v|^4\right)^{\frac{1}{2}}.
\end{align*}
Since $\tau\le1$, \eqref{eq-base-double-moment} with $2c$ and
\eqref{eq-gradient4-window} bound the product by $\sqrt{K_5c_4(2c,\eta)}$,
which gives \eqref{eq-weightedaverage} after enlarging $K_9(c,\eta)$ if necessary.
\end{proof}

For $\lambda>b$, the term $(\lambda-b)\int_{\Omega}u^2e^{\lambda v}$
in \eqref{eq-low-identity} is positive. To estimate this term, we apply the embedding
$W^{1,1}(\Omega)\hookrightarrow L^2(\Omega)$ to $ue^{\frac{\lambda v}{2}}$.
For the term containing $\nabla u$, we use H\"older's inequality to obtain
\[
\left(\int_{\Omega}e^{\frac{\lambda v}{2}}|\nabla u|\right)^2
\le\left(\int_{\Omega}\frac{|\nabla u|^2}{u}\right)\int_{\Omega}ue^{\lambda v}.
\]
The second factor on the right is part of the weighted energy. To apply
Lemma \ref{lem-uniform-gronwall}, we therefore need a uniform bound for
$\int_t^{t+\tau}\int_{\Omega}\frac{|\nabla u|^2}{u}$.
Under \eqref{eq-value-assumption} we estimate the entropy $\int_{\Omega}u\log u$
directly. Under \eqref{eq-main-assumption} we also use the second equation
of \eqref{eq-model} to cancel the term containing $\nabla u\cdot\nabla v$.

\begin{lemma}\label{lem-value-entropy}
Assume \eqref{eq-value-assumption}. Then there is a constant
$K_{10}>0$ such that
\begin{equation}\label{eq-value-entropy-bound}
\int_{\Omega}u\log u\le K_{10},
\qquad 0<t<T_{\max},
\end{equation}
and
\begin{equation}\label{eq-value-entropy-dissipation}
\int_t^{t+\tau}\int_{\Omega}\frac{|\nabla u|^2}{u}\le K_{10},
\qquad 0<t<\widetilde T_{\max}.
\end{equation}
\end{lemma}
\begin{proof}
Multiplying the first equation in \eqref{eq-model} by $\log u+1$,
integrating by parts and using the Neumann boundary conditions, we obtain
\begin{align}\label{eq-value-log-identity}
\frac{d}{dt}\int_{\Omega}u\log u
=&\int_{\Omega}(\log u+1)u_t\notag\\
=&\int_{\Omega}(\log u+1)
\{\Delta u+\nabla\cdot(u\chi(v)\nabla v)+au-bu^2\}\notag\\
=&-\int_{\Omega}\frac{|\nabla u|^2}{u}
-\int_{\Omega}\chi(v)\nabla u\cdot\nabla v
+\int_{\Omega}(au-bu^2)(\log u+1).
\end{align}
By Young's inequality and \eqref{eq-value-assumption}, there exists $c_1>0$ such that
\begin{align}\label{eq-value-entropy-cross}
\left|\int_{\Omega}\chi(v)\nabla u\cdot\nabla v\right|
\le&\frac{1}{2}\int_{\Omega}\frac{|\nabla u|^2}{u}
+\frac{1}{2}\int_{\Omega}u\chi(v)^2|\nabla v|^2\notag\\
\le&\frac{1}{2}\int_{\Omega}\frac{|\nabla u|^2}{u}
+\frac{c_1}{2}\int_{\Omega}u\chi(v)e^{bv}|\nabla v|^2.
\end{align}
There exists a constant $c_2>0$, depending only on $a,b$, such that
\[
(az-bz^2)(\log z+1)\le c_2,\qquad z>0.
\]
Substituting this bound and \eqref{eq-value-entropy-cross} into
\eqref{eq-value-log-identity}, we obtain
\begin{equation}\label{eq-value-entropy-inequality}
\frac{d}{dt}\int_{\Omega}u\log u
+\frac{1}{2}\int_{\Omega}\frac{|\nabla u|^2}{u}
\le\frac{c_1}{2}\int_{\Omega}u\chi(v)e^{bv}|\nabla v|^2+c_2|\Omega|.
\end{equation}
Taking $\lambda=b$ in \eqref{eq-low-dissipation} and using $\tau\le1$,
we find that the integral of the right-hand side over $(t,t+\tau)$ is bounded
uniformly in $t$. Also, $z\log z\ge-\frac{1}{e}$ for $z\ge0$, so
\[
\int_{\Omega}u\log u\ge-\frac{|\Omega|}{e}.
\]
Since $z\log z\le c_3(1+z^2)$ for $z\ge0$, \eqref{eq-u2-window} gives
\[
\int_t^{t+\tau}\int_{\Omega}u\log u\le c_3(|\Omega|+K_2),
\qquad 0<t<\widetilde T_{\max}.
\]
Applying Lemma \ref{lem-uniform-gronwall} to \eqref{eq-value-entropy-inequality}
with $y=\int_{\Omega}u\log u+\frac{|\Omega|}{e}$ and $h=0$ gives
\eqref{eq-value-entropy-bound} for $t>\tau$. On $[0,\tau]$ it follows
from the local boundedness of $u$.
It remains to integrate \eqref{eq-value-entropy-inequality} over $(t,t+\tau)$.
Using \eqref{eq-value-entropy-bound} at $t$, the lower bound $-\frac{|\Omega|}{e}$
at $t+\tau$ and \eqref{eq-low-dissipation} with $\lambda=b$, we obtain
\eqref{eq-value-entropy-dissipation} after enlarging $K_{10}$ if necessary.
\end{proof}

Under \eqref{eq-main-assumption}, we add $\frac{1}{2}\int_{\Omega}\chi(v)|\nabla v|^2$
to $\int_{\Omega}u\log u$ to cancel the mixed term involving $\chi(v)\nabla u\cdot\nabla v$.

\begin{lemma}\label{lem-entropy}
Assume \eqref{eq-main-assumption}. Then there is a constant
$K_{11}>0$ such that
\begin{equation}\label{eq-entropybound}
\int_{\Omega}u\log u+\frac{1}{2}\int_{\Omega}\chi(v)|\nabla v|^2
\le K_{11},\qquad 0<t<T_{\max},
\end{equation}
and
\begin{equation}\label{eq-entropy-dissipation}
\int_t^{t+\tau}\int_{\Omega}\frac{|\nabla u|^2}{u}\le K_{11},
\qquad 0<t<\widetilde T_{\max}.
\end{equation}
\end{lemma}
\begin{proof}
Equation \eqref{eq-value-log-identity} holds for every $\chi$
satisfying \eqref{eq-data}. Using the second equation in \eqref{eq-model} and
integrating by parts, we obtain
\begin{align}\label{eq-entropy-gradient-expansion}
\frac{1}{2}\frac{d}{dt}\int_{\Omega}\chi(v)|\nabla v|^2
=&\frac{1}{2}\int_{\Omega}\chi'(v)v_t|\nabla v|^2
+\int_{\Omega}\chi(v)\nabla v\cdot\nabla v_t\notag\\
=&-\int_{\Omega}\chi(v)v_t\Delta v
-\frac{1}{2}\int_{\Omega}\chi'(v)v_t|\nabla v|^2\notag\\
=&-\int_{\Omega}\chi(v)|\Delta v|^2
+\int_{\Omega}v\chi(v)\Delta v-\int_{\Omega}u\chi(v)\Delta v\notag\\
&-\frac{1}{2}\int_{\Omega}\chi'(v)|\nabla v|^2\Delta v
+\frac{1}{2}\int_{\Omega}v\chi'(v)|\nabla v|^2
-\frac{1}{2}\int_{\Omega}u\chi'(v)|\nabla v|^2.
\end{align}
Integrating by parts once more, we have
\[
\int_{\Omega}v\chi(v)\Delta v
=-\int_{\Omega}(\chi(v)+v\chi'(v))|\nabla v|^2
\]
and
\[
-\int_{\Omega}u\chi(v)\Delta v
=\int_{\Omega}\chi(v)\nabla u\cdot\nabla v
+\int_{\Omega}u\chi'(v)|\nabla v|^2.
\]
Substituting these identities into \eqref{eq-entropy-gradient-expansion}
gives
\begin{align}\label{eq-entropy-gradient-identity}
\frac{1}{2}\frac{d}{dt}\int_{\Omega}\chi(v)|\nabla v|^2
=&-\int_{\Omega}\chi(v)|\Delta v|^2
-\frac{1}{2}\int_{\Omega}\chi'(v)|\nabla v|^2\Delta v\notag\\
&-\int_{\Omega}\left(\chi(v)+\frac{v\chi'(v)}{2}\right)|\nabla v|^2
+\int_{\Omega}\chi(v)\nabla u\cdot\nabla v\notag\\
&+\frac{1}{2}\int_{\Omega}u\chi'(v)|\nabla v|^2.
\end{align}
When we add \eqref{eq-value-log-identity} and \eqref{eq-entropy-gradient-identity},
the terms involving $\chi(v)\nabla u\cdot\nabla v$ cancel, leaving
\begin{align}\label{eq-entropyid}
&\frac{d}{dt}\left(\int_{\Omega}u\log u
+\frac{1}{2}\int_{\Omega}\chi(v)|\nabla v|^2\right)
+\int_{\Omega}\frac{|\nabla u|^2}{u}
+\int_{\Omega}\chi(v)|\Delta v|^2\notag\\
=&-\frac{1}{2}\int_{\Omega}\chi'(v)|\nabla v|^2\Delta v
-\int_{\Omega}\left(\chi(v)+\frac{v\chi'(v)}{2}\right)|\nabla v|^2\notag\\
&+\frac{1}{2}\int_{\Omega}u\chi'(v)|\nabla v|^2
+\int_{\Omega}(au-bu^2)(\log u+1).
\end{align}
By \eqref{eq-main-assumption} and Young's inequality,
\begin{align*}
\frac{1}{2}\int_{\Omega}u|\chi'(v)||\nabla v|^2
\le&\frac{L}{2}\int_{\Omega}u\chi(v)e^{\frac{bv}{2}}|\nabla v|^2
+\frac{L}{2}\int_{\Omega}ue^{\frac{bv}{2}}|\nabla v|^2\\
\le&\frac{L}{2}\int_{\Omega}u\chi(v)e^{\frac{bv}{2}}|\nabla v|^2
+\frac{L}{4}\int_{\Omega}u^2
+\frac{L}{4}\int_{\Omega}e^{bv}|\nabla v|^4.
\end{align*}
After integration over $(t,t+\tau)$, we estimate the first and third terms by
\eqref{eq-low-dissipation} with $\lambda=\frac{b}{2}$ and the second by
\eqref{eq-u2-window}. For the term containing $v\chi'(v)$, we have
\[
-\frac{1}{2}\int_{\Omega}v\chi'(v)|\nabla v|^2
\le\frac{L}{2}\int_{\Omega}ve^{\frac{bv}{2}}(1+\chi(v))|\nabla v|^2.
\]
By \eqref{eq-automatic-growth}, there exist $c_1,c_2>0$ such that
$(1+s)e^{\frac{bs}{2}}(1+\chi(s))\le c_1\exp\{c_2e^{\frac{bs}{2}}\}$ for $s\ge0$, so \eqref{eq-weightedaverage} with $\eta=\frac{b}{2}$
bounds its integral over $(t,t+\tau)$.

Applying Young's inequality to the term containing $\Delta v$ and using
\eqref{eq-main-assumption} to bound $|\chi'(v)|$, we find
\[
\frac{1}{2}\left|\int_{\Omega}\chi'(v)|\nabla v|^2\Delta v\right|
\le\frac{1}{4}\int_{\Omega}(1+\chi(v))|\Delta v|^2
+c_3\int_{\Omega}(1+\chi(v))e^{bv}|\nabla v|^4.
\]
For $\rho(s)=1+\chi(s)$ and $R(s)=e^{\frac{bs}{2}}$,
\[
\frac{|\rho'(s)|}{\rho(s)}\le LR(s),\qquad
\frac{|R'(s)|}{R(s)}=\frac{b}{2}\le\frac{b}{2}R(s),\qquad
\int_{\Omega}R(v)^2|\nabla v|^2\le K_8(b)
\]
by \eqref{eq-main-assumption} and \eqref{eq-low-pointwise} with
$\lambda=b$. Applying Lemma \ref{lem-variable} with
$\varepsilon=\frac{1}{4c_3}$ therefore gives
\begin{align}\label{eq-entropy-mixed-bound}
\frac{1}{2}\left|\int_{\Omega}\chi'(v)|\nabla v|^2\Delta v\right|
\le&\frac{1}{2}\int_{\Omega}\chi(v)|\Delta v|^2
+\frac{1}{2}\int_{\Omega}|\Delta v|^2\notag\\
&+c_4\left(1+\int_{\Omega}e^{bv}|\Delta v|^2
+\int_{\Omega}e^{2bv}|\nabla v|^4\right)
\int_{\Omega}(1+\chi(v))|\nabla v|^2.
\end{align}
Since $v\ge0$, we have $|\Delta v|^2\le e^{bv}|\Delta v|^2$ and
$e^{bv}|\nabla v|^4\le e^{2bv}|\nabla v|^4$. Moreover, \eqref{eq-signal-H1} gives
\[
\int_{\Omega}(1+\chi(v))|\nabla v|^2
\le K_3+\int_{\Omega}\chi(v)|\nabla v|^2.
\]
The two remaining terms satisfy
\[
\int_{\Omega}(au-bu^2)(\log u+1)\le c_5,\qquad
-\int_{\Omega}\chi(v)|\nabla v|^2\le0.
\]
Since $u\log u\ge-\frac{1}{e}$ and $\chi\ge0$, we also have
\begin{equation}\label{eq-entropy}
\begin{gathered}
\int_{\Omega}u\log u+\frac{1}{2}\int_{\Omega}\chi(v)|\nabla v|^2
\ge-\frac{|\Omega|}{e},\\
\frac{1}{2}\int_{\Omega}\chi(v)|\nabla v|^2
\le\int_{\Omega}u\log u
+\frac{1}{2}\int_{\Omega}\chi(v)|\nabla v|^2+\frac{|\Omega|}{e}.
\end{gathered}
\end{equation}
In \eqref{eq-entropyid}, we use \eqref{eq-entropy-mixed-bound} to estimate the
term containing $\Delta v$. Substituting the bounds for the remaining terms,
using \eqref{eq-entropy} and absorbing $\frac{1}{2}\int_{\Omega}\chi(v)|\Delta v|^2$
on the left, we obtain
\begin{align}\label{eq-entropy-inequality}
&\frac{d}{dt}\left(\int_{\Omega}u\log u
+\frac{1}{2}\int_{\Omega}\chi(v)|\nabla v|^2\right)
+\int_{\Omega}\frac{|\nabla u|^2}{u}\notag\\
\le& c_6\left(1+\int_{\Omega}e^{bv}|\Delta v|^2
+\int_{\Omega}e^{2bv}|\nabla v|^4\right)
\left(\int_{\Omega}u\log u
+\frac{1}{2}\int_{\Omega}\chi(v)|\nabla v|^2
+\frac{|\Omega|}{e}\right)\notag\\
&+c_6\left(1+\int_{\Omega}e^{bv}|\Delta v|^2
+\int_{\Omega}e^{2bv}|\nabla v|^4\right)\notag\\
&+c_6\int_{\Omega}\left\{
u\chi(v)e^{\frac{bv}{2}}|\nabla v|^2+u^2
+ve^{\frac{bv}{2}}(1+\chi(v))|\nabla v|^2\right\}.
\end{align}
By \eqref{eq-low-dissipation} with $\lambda=b$, the coefficient of
$\int_{\Omega}u\log u+\frac{1}{2}\int_{\Omega}\chi(v)|\nabla v|^2+\frac{|\Omega|}{e}$
has a uniformly bounded integral over $(t,t+\tau)$. For the additive terms in
\eqref{eq-entropy-inequality}, we use \eqref{eq-low-dissipation} with
$\lambda=b$ and $\lambda=\frac{b}{2}$, \eqref{eq-u2-window} and
\eqref{eq-weightedaverage} with $\eta=\frac{b}{2}$ in view of \eqref{eq-automatic-growth}.
Since $u\log u\le c_7(1+u^2)$ and
$\chi(s)\le c_1\exp\{c_2e^{\frac{bs}{2}}\}$,
\eqref{eq-u2-window} and \eqref{eq-weightedaverage} with
$\eta=\frac{b}{2}$ give
\begin{equation}\label{eq-entropy-average}
\int_t^{t+\tau}\left(\int_{\Omega}u\log u
+\frac{1}{2}\int_{\Omega}\chi(v)|\nabla v|^2\right)\le c_8,
\qquad 0<t<\widetilde T_{\max}.
\end{equation}

After adding $\frac{|\Omega|}{e}$ to the entropy, we use \eqref{eq-entropy} and
\eqref{eq-entropy-average} to apply Lemma \ref{lem-uniform-gronwall} to
\eqref{eq-entropy-inequality} and obtain \eqref{eq-entropybound} for $t>\tau$.
On $[0,\tau]$ the bound follows from Lemma \ref{lem-local} and \eqref{eq-signal-H1}. Finally, we integrate
\eqref{eq-entropy-inequality} over $(t,t+\tau)$, using \eqref{eq-entropybound}
at $t$ and \eqref{eq-entropy} at $t+\tau$, to obtain \eqref{eq-entropy-dissipation}.
\end{proof}

Under either hypothesis of Theorem \ref{thm-main}, the entropy
dissipation has a uniformly bounded integral over $(t,t+\tau)$
by \eqref{eq-value-entropy-dissipation} or \eqref{eq-entropy-dissipation}.
We use this bound to estimate the positive term
$(q-b)\int_{\Omega}u^2e^{qv}$ when $q>b$.

\begin{lemma}\label{lem-all-exp}
Let $(u,v)$ be the solution given by Lemma \ref{lem-local}, and suppose that
$\chi$ satisfies \eqref{eq-value-assumption} or \eqref{eq-main-assumption}.
Then for every $q>0$ there exists $K_{12}(q)>0$ such that
\begin{align}\label{eq-all-expweights}
\int_{\Omega}ue^{qv}+\int_{\Omega}e^{qv}|\nabla v|^2
\le&K_{12}(q),\qquad\text{for all }t\in[0,T_{\max}),
\end{align}
and
\begin{align}\label{eq-all-exp-dissipation}
\int_t^{t+\tau}\int_{\Omega}e^{qv}|\Delta v|^2
+\int_t^{t+\tau}\int_{\Omega}e^{2qv}|\nabla v|^4
\le&K_{12}(q),\qquad\text{for all }t\in(0,\widetilde T_{\max}).
\end{align}
\end{lemma}
\begin{proof}
For $0<q\le b$ both estimates follow from Lemma \ref{lem-low}, with
the bound at $t=0$ given by the initial data. Let $q>b$.
Setting $\lambda=q$ in \eqref{eq-low-identity} and separating
the two terms in $(q-b)\int_{\Omega}u^2e^{qv}$, we obtain
\begin{align}\label{eq-q-energy-identity}
&\frac{d}{dt}\left(\int_{\Omega}ue^{qv}
+q\int_{\Omega}e^{qv}|\nabla v|^2\right)
+2q\int_{\Omega}e^{qv}|\Delta v|^2\notag\\
=&q\int_{\Omega}u^2e^{qv}
-q^2\int_{\Omega}e^{qv}|\nabla v|^2\Delta v\notag\\
&-q\int_{\Omega}(2+qv)e^{qv}|\nabla v|^2
-q\int_{\Omega}u\chi(v)e^{qv}|\nabla v|^2\notag\\
&-b\int_{\Omega}u^2e^{qv}+\int_{\Omega}u(a-qv)e^{qv}.
\end{align}
To estimate $q\int_{\Omega}u^2e^{qv}$, we apply the two-dimensional
Sobolev inequality to $ue^{\frac{qv}{2}}$. Using
\[
\nabla(ue^{\frac{qv}{2}})=e^{\frac{qv}{2}}\nabla u+\frac{q}{2}ue^{\frac{qv}{2}}\nabla v,
\]
we obtain, for some $c_1>0$,
\begin{align*}
\int_{\Omega}u^2e^{qv}
\le& 2c_1\left(\int_{\Omega}e^{\frac{qv}{2}}|\nabla u|\right)^2
+\frac{q^2c_1}{2}\left(\int_{\Omega}ue^{\frac{qv}{2}}|\nabla v|\right)^2
+c_1\left(\int_{\Omega}ue^{\frac{qv}{2}}\right)^2.
\end{align*}
By H\"older's inequality and \eqref{eq-mass},
\begin{align*}
\left(\int_{\Omega}e^{\frac{qv}{2}}|\nabla u|\right)^2
\le&\left(\int_{\Omega}\frac{|\nabla u|^2}{u}\right)
       \int_{\Omega}ue^{qv},\\
\left(\int_{\Omega}ue^{\frac{qv}{2}}|\nabla v|\right)^2
\le&\left(\int_{\Omega}u\right)
       \int_{\Omega}ue^{qv}|\nabla v|^2
\le K_1\int_{\Omega}ue^{qv}|\nabla v|^2,\\
\left(\int_{\Omega}ue^{\frac{qv}{2}}\right)^2
\le&\left(\int_{\Omega}u\right)\int_{\Omega}ue^{qv}
\le K_1\int_{\Omega}ue^{qv}.
\end{align*}
Consequently,
\[
\int_{\Omega}u^2e^{qv}
\le c_2\left(1+\int_{\Omega}\frac{|\nabla u|^2}{u}\right)
       \int_{\Omega}ue^{qv}
+c_3(q)\int_{\Omega}ue^{qv}|\nabla v|^2.
\]
Applying Young's inequality to the last term, we obtain
\[
q\int_{\Omega}u^2e^{qv}
\le 2qc_2\left(1+\int_{\Omega}\frac{|\nabla u|^2}{u}\right)\int_{\Omega}ue^{qv}
+c_4(q)\int_{\Omega}e^{qv}|\nabla v|^4.
\]
In view of \eqref{eq-signal-H1}, Lemma \ref{lem-variable} applies with
$\rho(s)=e^{qs}$, $R=1$ and $\varepsilon=\frac{q}{2c_4(q)}$, giving
\begin{align*}
c_4(q)\int_{\Omega}e^{qv}|\nabla v|^4
\le&\frac{q}{2}\int_{\Omega}e^{qv}|\Delta v|^2\\
&+c_5(q)\left(1+\int_{\Omega}|\Delta v|^2
+\int_{\Omega}|\nabla v|^4\right)\int_{\Omega}e^{qv}|\nabla v|^2.
\end{align*}
Combining the last two estimates, we find
\begin{align}\label{eq-ordinary-source}
&q\int_{\Omega}u^2e^{qv}\notag\\
\le&2qc_2\left(1+\int_{\Omega}\frac{|\nabla u|^2}{u}\right)\int_{\Omega}ue^{qv}
+\frac{q}{2}\int_{\Omega}e^{qv}|\Delta v|^2\notag\\
&+c_5(q)\left(1+\int_{\Omega}|\Delta v|^2+\int_{\Omega}|\nabla v|^4\right)
\int_{\Omega}e^{qv}|\nabla v|^2\notag\\
\le&\frac{q}{2}\int_{\Omega}e^{qv}|\Delta v|^2\notag\\
&+c_6(q)
\left(1+\int_{\Omega}|\Delta v|^2+\int_{\Omega}|\nabla v|^4
          +\int_{\Omega}\frac{|\nabla u|^2}{u}\right)
\left(\int_{\Omega}ue^{qv}+q\int_{\Omega}e^{qv}|\nabla v|^2\right).
\end{align}
We next apply Young's inequality to the term containing $|\nabla v|^2\Delta v$
in \eqref{eq-q-energy-identity} and then use Lemma \ref{lem-variable} with
$\rho(s)=e^{qs}$, $R=1$ and $\varepsilon=\frac{1}{4q^2}$ to obtain
\begin{align}\label{eq-exponential-mixed-bound}
&q^2\left|\int_{\Omega}e^{qv}|\nabla v|^2\Delta v\right|\notag\\
\le&\frac{q}{4}\int_{\Omega}e^{qv}|\Delta v|^2
+q^3\int_{\Omega}e^{qv}|\nabla v|^4\notag\\
\le&\frac{q}{2}\int_{\Omega}e^{qv}|\Delta v|^2
+c_7(q)\left(1+\int_{\Omega}|\Delta v|^2
                    +\int_{\Omega}|\nabla v|^4\right)
                  \int_{\Omega}e^{qv}|\nabla v|^2.
\end{align}
Since $u,v,\chi(v)\ge0$, the remaining terms satisfy
\begin{align*}
-q\int_{\Omega}(2+qv)e^{qv}|\nabla v|^2\le&0,\\
-q\int_{\Omega}u\chi(v)e^{qv}|\nabla v|^2\le&0,\\
-b\int_{\Omega}u^2e^{qv}\le&0,\\
\int_{\Omega}u(a-qv)e^{qv}\le& a\int_{\Omega}ue^{qv}.
\end{align*}
Substituting these four bounds into \eqref{eq-q-energy-identity} and using
\eqref{eq-ordinary-source} and \eqref{eq-exponential-mixed-bound}, we obtain
\begin{align*}
&\frac{d}{dt}\left(\int_{\Omega}ue^{qv}+q\int_{\Omega}e^{qv}|\nabla v|^2\right)
+2q\int_{\Omega}e^{qv}|\Delta v|^2\\
\le&q\int_{\Omega}u^2e^{qv}
+q^2\left|\int_{\Omega}e^{qv}|\nabla v|^2\Delta v\right|
+a\int_{\Omega}ue^{qv}\\
\le&\frac{q}{2}\int_{\Omega}e^{qv}|\Delta v|^2
+c_6(q)\left(1+\int_{\Omega}|\Delta v|^2
+\int_{\Omega}|\nabla v|^4+\int_{\Omega}\frac{|\nabla u|^2}{u}\right)
\left(\int_{\Omega}ue^{qv}+q\int_{\Omega}e^{qv}|\nabla v|^2\right)\\
&+\frac{q}{2}\int_{\Omega}e^{qv}|\Delta v|^2
+c_7(q)\left(1+\int_{\Omega}|\Delta v|^2+\int_{\Omega}|\nabla v|^4\right)
\int_{\Omega}e^{qv}|\nabla v|^2+a\int_{\Omega}ue^{qv}.
\end{align*}
We absorb the two integrals containing $|\Delta v|^2$ into the left-hand side. With
$c_8(q)=c_6(q)+\frac{c_7(q)}{q}+a$, the resulting differential inequality is
\begin{align}\label{eq-ordinary-frozen}
&\frac{d}{dt}\left(\int_{\Omega}ue^{qv}
                    +q\int_{\Omega}e^{qv}|\nabla v|^2\right)\notag\\
&+q\int_{\Omega}e^{qv}|\Delta v|^2\notag\\
\le& c_8(q)\left(1+\int_{\Omega}|\Delta v|^2+\int_{\Omega}|\nabla v|^4
                  +\int_{\Omega}\frac{|\nabla u|^2}{u}\right)
 \left(\int_{\Omega}ue^{qv}
                   +q\int_{\Omega}e^{qv}|\nabla v|^2\right).
\end{align}
By \eqref{eq-expaverage} with $\lambda=q$, the energy in
\eqref{eq-ordinary-frozen} has a bounded integral over $(t,t+\tau)$.
Under \eqref{eq-value-assumption}, \eqref{eq-laplace-window},
\eqref{eq-gradient4-window} and \eqref{eq-value-entropy-dissipation}
bound the integral of its coefficient over the same interval by
$c_8(q)(\tau+K_4+K_5+K_{10})$. Under \eqref{eq-main-assumption},
\eqref{eq-entropy-dissipation} gives the same bound with $K_{11}$ in place of $K_{10}$.
Applying Lemma \ref{lem-uniform-gronwall} for $t>\tau$ and arguing on
$[0,\tau]$ as in the proof of Lemma \ref{lem-low}, we obtain
\eqref{eq-all-expweights} with upper bound $c_9(q)$.

Integrating \eqref{eq-ordinary-frozen} over $(t,t+\tau)$ and using
\eqref{eq-all-expweights} and the integral bound for the coefficient, we obtain
\begin{equation}\label{eq-exponential-laplacian-window}
\int_t^{t+\tau}\int_{\Omega}e^{qv}|\Delta v|^2\le c_{10}(q).
\end{equation}
Using \eqref{eq-all-expweights} in place of \eqref{eq-low-pointwise},
we repeat the interpolation in
\eqref{eq-exponential-quartic} with $\lambda=q$. Integrating over $(t,t+\tau)$, we obtain
\[
\int_t^{t+\tau}\int_{\Omega}e^{2qv}|\nabla v|^4
\le c_{11}(q)\left(\tau+\int_t^{t+\tau}\int_{\Omega}e^{qv}|\Delta v|^2
+\int_t^{t+\tau}\int_{\Omega}|\nabla v|^4\right)\le c_{12}(q),
\]
by \eqref{eq-exponential-laplacian-window} and \eqref{eq-gradient4-window}.
Choosing $K_{12}(q)=\max\{c_9(q),c_{10}(q)+c_{12}(q)\}$ proves
\eqref{eq-all-expweights} and \eqref{eq-all-exp-dissipation}.
\end{proof}

\section{Estimates with double-exponential weights}\label{sec-double-exponential}
Under \eqref{eq-value-assumption}, taking $q=2b$ in
\eqref{eq-all-exp-dissipation} suffices to control
$\int_t^{t+\tau}\int_{\Omega}\chi(v)^4|\nabla v|^4$.
Throughout this section we assume \eqref{eq-main-assumption}. In view of
\eqref{eq-automatic-growth}, we consider the weight
\begin{equation}\label{eq-weight-family}
w(s)=\exp\{Ae^{\frac{bs}{2}}\},\qquad s\ge0,\quad A\ge1.
\end{equation}
For sufficiently large $A$, we have $\chi(s)^4\le(1+\chi(0))^4w(s)$ for all $s\ge0$,
so it suffices to bound $\int_t^{t+\tau}\int_{\Omega}w(v)|\nabla v|^4$ uniformly in $t$.
As in Section \ref{sec-exponential}, we begin with a differential inequality for the weighted energy
$\int_{\Omega}uw(v)+\int_{\Omega}w'(v)|\nabla v|^2$.

\begin{lemma}\label{lem-high-energy}
There exist a constant $K_{13}>0$ and a positive nondecreasing function
$K_{14}(A)$ on $[1,\infty)$ such that, for every $A\ge1$,
\begin{align}\label{eq-fixed-weight-energy}
&\frac{d}{dt}\left(\int_{\Omega}uw(v)  +\int_{\Omega}w'(v)|\nabla v|^2\right)+\int_{\Omega}w'(v)|\Delta v|^2\notag\\
\le& \left[a+K_{14}(A)
 \left(1+\int_{\Omega}e^{bv}|\Delta v|^2
          +\int_{\Omega}e^{2bv}|\nabla v|^4\right)\right]
 \left(\int_{\Omega}uw(v)
                  +\int_{\Omega}w'(v)|\nabla v|^2\right)\notag\\
&+K_{13}\left(1+\int_{\Omega}\frac{|\nabla u|^2}{u}\right)
                 \int_{\Omega}uw'(v),
\end{align}
for all $t\in(0,T_{\max})$.
\end{lemma}
\begin{proof}
Fix $A\ge1$. Using \eqref{eq-model} and integrating by parts under the
Neumann boundary conditions, we obtain
\begin{align}
\frac{d}{dt}\int_{\Omega}uw(v)
=&\int_{\Omega}w(v)
       [\Delta u+\nabla\cdot(u\chi(v)\nabla v)+au-bu^2]\notag\\
&+\int_{\Omega}uw'(v)(\Delta v-v+u)\notag\\
=&-\int_{\Omega}w'(v)\nabla u\cdot\nabla v
-\int_{\Omega}u\chi(v)w'(v)|\nabla v|^2\notag\\
&-\int_{\Omega}\nabla(uw'(v))\cdot\nabla v
+\int_{\Omega}(au-bu^2)w(v)\notag\\
&-\int_{\Omega}uvw'(v)+\int_{\Omega}u^2w'(v)\notag\\
=&-2\int_{\Omega}w'(v)\nabla u\cdot\nabla v
-\int_{\Omega}u(w''(v)+\chi(v)w'(v))
                          |\nabla v|^2\notag\\
&+\int_{\Omega}u^2(w'(v)-bw(v))
+\int_{\Omega}u(aw(v)-vw'(v)).
\label{eq-uw-identity}
\end{align}
To cancel the term containing $\nabla u\cdot\nabla v$, we next
differentiate $\int_{\Omega}w'(v)|\nabla v|^2$. By the second equation in
\eqref{eq-model} and integration by parts, we have
\begin{align}\label{eq-weight-gradient-expansion}
\frac{d}{dt}\int_{\Omega}w'(v)|\nabla v|^2
=&\int_{\Omega}w''(v)v_t|\nabla v|^2
+2\int_{\Omega}w'(v)\nabla v\cdot\nabla v_t\notag\\
=&-2\int_{\Omega}w'(v)(\Delta v)v_t
  -\int_{\Omega}w''(v)|\nabla v|^2v_t\notag\\
=&-2\int_{\Omega}w'(v)\Delta v(\Delta v-v+u)\notag\\
&-\int_{\Omega}w''(v)|\nabla v|^2(\Delta v-v+u).
\end{align}
We integrate the terms containing $v\Delta v$ and $u\Delta v$ by parts once more.
\begin{align*}
2\int_{\Omega}vw'(v)\Delta v
=&-2\int_{\Omega}(w'(v)+vw''(v))|\nabla v|^2,\\
-2\int_{\Omega}uw'(v)\Delta v
=&2\int_{\Omega}w'(v)\nabla u\cdot\nabla v
+2\int_{\Omega}uw''(v)|\nabla v|^2.
\end{align*}
Substituting these identities into \eqref{eq-weight-gradient-expansion}, we obtain
\begin{align}\label{eq-gradient-identity}
\frac{d}{dt}\int_{\Omega}w'(v)|\nabla v|^2
=&-2\int_{\Omega}w'(v)|\Delta v|^2
-\int_{\Omega}w''(v)|\nabla v|^2\Delta v\notag\\
&-\int_{\Omega}(2w'(v)+vw''(v))|\nabla v|^2\notag\\
&+2\int_{\Omega}w'(v)\nabla u\cdot\nabla v
+\int_{\Omega}uw''(v)|\nabla v|^2.
\end{align}
Adding \eqref{eq-uw-identity} and \eqref{eq-gradient-identity} and cancelling the
terms containing $\nabla u\cdot\nabla v$ and $uw''(v)|\nabla v|^2$, we obtain
\begin{align}\label{eq-mixed-high-identity}
&\frac{d}{dt}\left(\int_{\Omega}uw(v)
                    +\int_{\Omega}w'(v)|\nabla v|^2\right)\notag\\
=&-2\int_{\Omega}w'(v)|\Delta v|^2
-\int_{\Omega}w''(v)|\nabla v|^2\Delta v\notag\\
&-\int_{\Omega}(2w'(v)+vw''(v))|\nabla v|^2
-\int_{\Omega}u\chi(v)w'(v)|\nabla v|^2\notag\\
&+\int_{\Omega}u^2w'(v)-b\int_{\Omega}u^2w(v)
+\int_{\Omega}u(aw(v)-vw'(v)).
\end{align}
We first estimate the positive term $\int_{\Omega}u^2w'(v)$ in
\eqref{eq-mixed-high-identity}. By \eqref{eq-weight-family},
\begin{equation}\label{eq-weight-derivatives}
w'(s)=\frac{Ab}{2}e^{\frac{bs}{2}}w(s),\qquad
w''(s)=\frac{b}{2}\bigl(1+Ae^{\frac{bs}{2}}\bigr)w'(s)
\le c_1(A)e^{\frac{bs}{2}}w'(s),
\end{equation}
where $c_1(A)=\frac{b}{2}(1+A)$. To apply the two-dimensional Sobolev
inequality to $u\sqrt{w'(v)}$, we write
\[
\nabla(u\sqrt{w'(v)})=\sqrt{w'(v)}\nabla u
+\frac{uw''(v)}{2\sqrt{w'(v)}}\nabla v.
\]
By H\"older's inequality and \eqref{eq-mass},
\begin{align*}
\left(\int_{\Omega}\sqrt{w'(v)}|\nabla u|\right)^2
\le&\left(\int_{\Omega}\frac{|\nabla u|^2}{u}\right)\int_{\Omega}uw'(v),\\
\left(\int_{\Omega}u\sqrt{w'(v)}\right)^2
\le& K_1\int_{\Omega}uw'(v).
\end{align*}
For the term containing $w''$, we combine H\"older's inequality with
\eqref{eq-weight-derivatives} and the bound for $\int_{\Omega}ue^{bv}$ in
\eqref{eq-low-pointwise} with $\lambda=b$ to obtain
\begin{align*}
\left(\int_{\Omega}u\frac{w''(v)}{\sqrt{w'(v)}}|\nabla v|\right)^2
\le&\left(\int_{\Omega}ue^{bv}\right)
\int_{\Omega}ue^{-bv}\frac{w''(v)^2}{w'(v)}|\nabla v|^2\\
\le& K_8(b)c_1(A)^2\int_{\Omega}uw'(v)|\nabla v|^2.
\end{align*}
Substituting these estimates into the Sobolev inequality
$W^{1,1}(\Omega)\hookrightarrow L^2(\Omega)$ gives
\[
\int_{\Omega}u^2w'(v)
\le c_2\left(1+\int_{\Omega}\frac{|\nabla u|^2}{u}\right)
        \int_{\Omega}uw'(v)
+c_3(A)\int_{\Omega}uw'(v)|\nabla v|^2,
\]
where $c_2$ depends only on $\Omega$ and $K_1$. By Young's inequality applied to the last
term, we obtain
\begin{equation}\label{eq-double-source-bound}
\int_{\Omega}u^2w'(v)
\le K_{13}\left(1+\int_{\Omega}\frac{|\nabla u|^2}{u}\right)
           \int_{\Omega}uw'(v)
+c_3(A)^2\int_{\Omega}w'(v)|\nabla v|^4,
\end{equation}
with $K_{13}=2c_2$ independent of $A$.

Applying Young's inequality and using \eqref{eq-weight-derivatives} to bound
$w''$, we find
\begin{equation}\label{eq-double-mixed-bound}
\left|\int_{\Omega}w''(v)|\nabla v|^2\Delta v\right|
\le\frac{1}{2}\int_{\Omega}w'(v)|\Delta v|^2
+\frac{c_1(A)^2}{2}\int_{\Omega}w'(v)e^{bv}|\nabla v|^4.
\end{equation}
Since $e^{bv}\ge1$, it remains in both cases to estimate
$\int_{\Omega}w'(v)e^{bv}|\nabla v|^4$. We apply Lemma \ref{lem-variable} with
$\rho(s)=w'(s)$ and $R(s)=e^{\frac{bs}{2}}$. Indeed,
\[
|\rho'(s)|\le c_1(A)R(s)\rho(s),\qquad
R'(s)=\frac{b}{2}R(s)\le\frac{b}{2}R(s)^2,\qquad
\int_{\Omega}R(v)^2|\nabla v|^2\le K_8(b),
\]
by \eqref{eq-weight-derivatives} and \eqref{eq-low-pointwise}.
Taking $\varepsilon=\frac{1}{2c_3(A)^2+c_1(A)^2}$ in \eqref{eq-variable}
and adding \eqref{eq-double-source-bound} and \eqref{eq-double-mixed-bound}, we obtain
\begin{align}\label{eq-source}
&\int_{\Omega}u^2w'(v)
+\left|\int_{\Omega}w''(v)|\nabla v|^2\Delta v\right|\notag\\
\le&\int_{\Omega}w'(v)|\Delta v|^2
+K_{13}\left(1+\int_{\Omega}\frac{|\nabla u|^2}{u}\right)
             \int_{\Omega}uw'(v)\notag\\
&+K_{14}(A)\left(1+\int_{\Omega}e^{bv}|\Delta v|^2
                   +\int_{\Omega}e^{2bv}|\nabla v|^4\right)
             \int_{\Omega}w'(v)|\nabla v|^2.
\end{align}

The same choices of $c_1(A)$, $c_3(A)$ and $\varepsilon$ work for every weight
with parameter in $[1,A]$, so we choose $K_{14}(A)$ in \eqref{eq-source}
to be nondecreasing in $A$.
Since $u,v,\chi(v)\ge0$ and $w,w',w''>0$,
the remaining terms satisfy
\begin{align*}
-\int_{\Omega}(2w'(v)+vw''(v))|\nabla v|^2\le&0,\\
-\int_{\Omega}u\chi(v)w'(v)|\nabla v|^2\le&0,\\
-b\int_{\Omega}u^2w(v)\le&0,\\
\int_{\Omega}u(aw(v)-vw'(v))\le& a\int_{\Omega}uw(v).
\end{align*}
Substituting these bounds and \eqref{eq-source} into
\eqref{eq-mixed-high-identity} gives \eqref{eq-fixed-weight-energy}.
\end{proof}

The last term in \eqref{eq-fixed-weight-energy} contains
$\int_{\Omega}uw'(v)$, whereas the energy contains $\int_{\Omega}uw(v)$.
The coefficient of $\int_{\Omega}uw'(v)$ is independent of $A$. The weight $w(v+h)$ has the same form as $w(v)$, with
$Ae^{\frac{bh}{2}}$ in place of $A$. We therefore let $h$ depend on time
and choose $h'$ to cancel the term containing $uw'(v+h)$.

\begin{lemma}
\label{lem-direct}
For the fixed $A\ge1$ in \eqref{eq-weight-family}, there is
$K_{15}(A)>0$ such that
\begin{align}\label{eq-direct-weighted}
\int_{\Omega}uw(v)+\int_{\Omega}w'(v)|\nabla v|^2
\le&K_{15}(A),\qquad\text{for all }t\in[0,T_{\max}),
\end{align}
and
\begin{align}\label{eq-direct-dissipation}
\int_t^{t+\tau}\int_{\Omega}w'(v)|\Delta v|^2
+\int_t^{t+\tau}\int_{\Omega}w(v)|\nabla v|^4
\le&K_{15}(A),\qquad\text{for all }t\in(0,\widetilde T_{\max}).
\end{align}
\end{lemma}
\begin{proof}
Fix $t\in(0,\widetilde T_{\max})$. We estimate the energy at times
$r\in[t,t+\tau]$. For $r>\tau$ we shall apply the uniform Gr\"onwall
lemma on $(r-\tau,r)$, so we define the shift on the larger interval
$(\max\{0,t-\tau\},t+\tau]$.

With $K_{13}$ from \eqref{eq-fixed-weight-energy}, set
\[
h(r)=K_{13}\int_r^{t+\tau}\left(1+\int_{\Omega}\frac{|\nabla u(x,s)|^2}{u(x,s)}dx\right)ds,
\qquad \max\{0,t-\tau\}<r\le t+\tau.
\]
If $r\ge t$, the interval $(r,t+\tau)$ is contained in $(t,t+\tau)$.
If $r<t$, the two intervals $(r,r+\tau)$ and $(t,t+\tau)$ cover
$(r,t+\tau)$. Applying \eqref{eq-entropy-dissipation} on these intervals
gives
\[
0\le h(r)\le K_{13}(2\tau+2K_{11})\le c_1,
\qquad h(t+\tau)=0.
\]
Here $c_1=2K_{13}(1+K_{11})$ is independent of $t$ and $A$.
By the definition of $h$, its derivative satisfies
\begin{equation}\label{eq-shift-derivative}
h'(r)=-K_{13}\left(1+\int_{\Omega}\frac{|\nabla u(x,r)|^2}{u(x,r)}dx\right),
\quad\text{for a.e. }r.
\end{equation}
At each fixed time $r$, the shifted weight satisfies
\[
w(s+h(r))=\exp\{(Ae^{\frac{bh(r)}{2}})e^{\frac{bs}{2}}\},\qquad s\ge0.
\]
Since the parameter $Ae^{\frac{bh(r)}{2}}$ lies between $A$ and
$Ae^{\frac{bc_1}{2}}$, Lemma \ref{lem-high-energy} applies with the same $K_{13}$ at each $r$.
Since $c_1$ is independent of $t$ and $A$ and $K_{14}$ is nondecreasing,
we use $c_2(A)=K_{14}(Ae^{\frac{bc_1}{2}})$ throughout the interval.

Since $0\le h\le c_1$ and $w,w'$ are increasing, replacing $h(r)$ by $c_1$ gives
an upper bound for the time integral of the energy. We apply Young's inequality to $uw(v+c_1)$
and note that for all $s\ge0$,
\[
w'(s+c_1)=\frac{Ab}{2}e^{\frac{bc_1}{2}}e^{\frac{bs}{2}}
                 \exp\{Ae^{\frac{bc_1}{2}}e^{\frac{bs}{2}}\}
\le c_3(A)\exp\{(Ae^{\frac{bc_1}{2}}+1)e^{\frac{bs}{2}}\}.
\]
Using \eqref{eq-u2-window} and \eqref{eq-base-double-moment}, we bound the integrals
of $u^2$ and $w(v+c_1)^2$ over $\Omega\times(s,s+\tau)$. For the integral
containing $|\nabla v|^2$, we apply \eqref{eq-weightedaverage}, again with
$\eta=\frac{b}{2}<b$. Since $\tau\le1$, we obtain
\begin{equation}\label{eq-fixed-energy-average}
\begin{aligned}
&\int_s^{s+\tau}\biggl[
 \int_{\Omega}u(x,r)w(v(x,r)+c_1)dx\\
&+\int_{\Omega}w'(v(x,r)+c_1)|\nabla v(x,r)|^2dx
 \biggr]dr\le c_4(A),
 \quad(0<s<\widetilde T_{\max}).
\end{aligned}
\end{equation}
Likewise, \eqref{eq-low-dissipation} with $\lambda=b$ gives
\begin{align}\label{eq-shift-coefficient-average}
&\int_s^{s+\tau}\biggl[1+\int_{\Omega}e^{bv(x,r)}|\Delta v(x,r)|^2dx\notag\\
&+\int_{\Omega}e^{2bv(x,r)}|\nabla v(x,r)|^4dx
 \biggr]dr\le 1+K_8(b),
\end{align}
for all $0<s<\widetilde T_{\max}$.

For this fixed $t$, set
\[
\begin{aligned}
I_1(r)=&\int_{\Omega}u(x,r)w(v(x,r)+h(r))dx\\
&+\int_{\Omega}w'(v(x,r)+h(r))|\nabla v(x,r)|^2dx.
\end{aligned}
\]
We differentiate $I_1(r)$ with respect to $r$, taking into account the dependence of the weight on $h(r)$.
\begin{align*}
I_1'(r)
=&\int_{\Omega}\partial_ru(x,r)w(v(x,r)+h(r))dx\\
&+\int_{\Omega}u(x,r)w'(v(x,r)+h(r))\partial_rv(x,r)dx\\
&+\int_{\Omega}w''(v(x,r)+h(r))
                   \partial_r v(x,r)|\nabla v(x,r)|^2dx\\
&+2\int_{\Omega}w'(v(x,r)+h(r))
                   \nabla v(x,r)\cdot\nabla\partial_r v(x,r)dx\\
&+h'(r)\int_{\Omega}u(x,r)w'(v(x,r)+h(r))dx\\
&+h'(r)\int_{\Omega}w''(v(x,r)+h(r))|\nabla v(x,r)|^2dx.
\end{align*}
The first four integrals give the derivative with the shift held fixed.
Applying \eqref{eq-fixed-weight-energy} to these terms with the weight
$s\mapsto w(s+h(r))$ and adding the two terms containing $h'(r)$, we obtain
\begin{align*}
I_1'(r)
&+\int_{\Omega}w'(v(x,r)+h(r))|\Delta v(x,r)|^2dx\\
\le& aI_1(r)
 +c_2(A)\biggl[1+\int_{\Omega}e^{bv(x,r)}|\Delta v(x,r)|^2dx\\
&+\int_{\Omega}e^{2bv(x,r)}|\nabla v(x,r)|^4dx
       \biggr]I_1(r)\\
&+\left[K_{13}\left(1+\int_{\Omega}\frac{|\nabla u(x,r)|^2}{u(x,r)}dx\right)
                 +h'(r)\right]
       \int_{\Omega}u(x,r)w'(v(x,r)+h(r))dx\\
&+h'(r)\int_{\Omega}w''(v(x,r)+h(r))|\nabla v(x,r)|^2dx.
\end{align*}
By \eqref{eq-shift-derivative}, the coefficient of
$\int_{\Omega}u(x,r)w'(v(x,r)+h(r))dx$ vanishes. Since $h'\le0$ and $w''\ge0$,
we drop the last term to obtain
\begin{align}\label{eq-closed-energy}
I_1'(r)
&+\int_{\Omega}w'(v(x,r)+h(r))|\Delta v(x,r)|^2dx\notag\\
\le& aI_1(r)
 +c_2(A)\biggl[1+\int_{\Omega}e^{bv(x,r)}|\Delta v(x,r)|^2dx\notag\\
&+\int_{\Omega}e^{2bv(x,r)}|\nabla v(x,r)|^4dx
       \biggr]I_1(r).
\end{align}

For $t<r<t+\tau$ with $r>\tau$, we use \eqref{eq-fixed-energy-average} and
\eqref{eq-shift-coefficient-average} to apply Lemma \ref{lem-uniform-gronwall}
to \eqref{eq-closed-energy} on $(r-\tau,r)$.
For $r\le\tau$, we bound $I_1(r)$ as in Lemma \ref{lem-low} using $0\le h\le c_1$.
Thus, by continuity and the monotonicity of $w$ and $w'$,
\begin{align*}
&\int_{\Omega}u(x,r)w(v(x,r))dx
 +\int_{\Omega}w'(v(x,r))|\nabla v(x,r)|^2dx\\
\le& I_1(r)\le c_5(A),\qquad t\le r\le t+\tau,
\end{align*}
where $c_5(A)$ is independent of $t$. Since $t$ was arbitrary and the
bound at time zero follows from the initial data, \eqref{eq-direct-weighted} holds.

Since $h\ge0$ and $I_1(t+\tau)\ge0$, integrating \eqref{eq-closed-energy}
over $(t,t+\tau)$ and using \eqref{eq-shift-coefficient-average} gives
\begin{equation}\label{eq-weighted-laplacian-window}
\int_t^{t+\tau}\int_{\Omega}w'(v(x,r))|\Delta v(x,r)|^2dxdr
\le I_1(t)+c_5(A)\bigl[a\tau+c_2(A)(1+K_8(b))\bigr]\le c_6(A).
\end{equation}
Finally, we apply Lemma \ref{lem-variable} with
$\rho(s)=w(s)$ and $R(s)=e^{\frac{bs}{2}}$. By
\eqref{eq-weight-derivatives}, $w'(s)=\frac{Ab}{2}R(s)w(s)$,
and \eqref{eq-low-pointwise} with $\lambda=b$ bounds
$\int_{\Omega}R(v)^2|\nabla v|^2$. Taking $\varepsilon=1$ in
\eqref{eq-variable} and using $R(v)\ge1$, we obtain
\begin{align}\label{eq-quartic-interpolation}
&\int_{\Omega}w(v(x,t))|\nabla v(x,t)|^4dx\notag\\*
\le&\int_{\Omega}w(v(x,t))|\Delta v(x,t)|^2dx\notag\\*
&+c_7(A)\biggl[1+\int_{\Omega}e^{bv(x,t)}|\Delta v(x,t)|^2dx\notag\\*
&+\int_{\Omega}e^{2bv(x,t)}|\nabla v(x,t)|^4dx
       \biggr]\int_{\Omega}w(v(x,t))|\nabla v(x,t)|^2dx.
\end{align}
We integrate \eqref{eq-quartic-interpolation} over $(t,t+\tau)$. Since
$w'\ge\frac{Ab}{2}w$, \eqref{eq-weighted-laplacian-window} bounds the time integral
of the first term, while \eqref{eq-direct-weighted} bounds
$\int_{\Omega}w(v)|\nabla v|^2$ uniformly in time. Using
\eqref{eq-shift-coefficient-average} for the remaining factor, we obtain
\[
\int_t^{t+\tau}\int_{\Omega}w(v(x,s))|\nabla v(x,s)|^4dxds\le c_8(A),
\qquad(0<t<\widetilde T_{\max}).
\]
The two integrals in \eqref{eq-direct-dissipation} are thus bounded by $c_6(A)$ and $c_8(A)$, respectively.
\end{proof}

\section{Proof of the main theorem}\label{sec-boundedness}
We now use the weighted estimates from Sections \ref{sec-exponential} and
\ref{sec-double-exponential} to obtain the following space-time $L^4$ bound for $\chi(v)\nabla v$.

\begin{lemma}\label{lem-fourth-power}
Under either hypothesis on $\chi$ in Theorem \ref{thm-main},
there exists $K_{16}>0$ such that
\begin{equation}\label{eq-fourth-power}
\int_t^{t+\tau}\int_{\Omega}\chi(v)^4|\nabla v|^4\le K_{16},
\qquad\text{for all }t\in(0,\widetilde T_{\max}).
\end{equation}
\end{lemma}
\begin{proof}
Suppose first that \eqref{eq-value-assumption} holds. Then
$\chi(v)^4\le C^4e^{4bv}$, and \eqref{eq-all-exp-dissipation} with $q=2b$ gives
\[
\int_t^{t+\tau}\int_{\Omega}\chi(v)^4|\nabla v|^4
\le C^4\int_t^{t+\tau}\int_{\Omega}e^{4bv}|\nabla v|^4\le C^4K_{12}(2b).
\]
Now suppose that \eqref{eq-main-assumption} holds and take
$A=\max\{1,\frac{8L}{b}\}$ in \eqref{eq-weight-family}. Using the bound
$\chi(s)^4\le(1+\chi(0))^4w(s)$ for $s\ge0$ from \eqref{eq-automatic-growth}
together with \eqref{eq-direct-dissipation}, we find
\[
\int_t^{t+\tau}\int_{\Omega}\chi(v)^4|\nabla v|^4
\le (1+\chi(0))^4\int_t^{t+\tau}\int_{\Omega}w(v)|\nabla v|^4
\le (1+\chi(0))^4K_{15}(A).
\]
This proves \eqref{eq-fourth-power} in both cases.
\end{proof}

We use Lemma \ref{lem-fourth-power} to bound $u$ in $L^2(\Omega)$.

\begin{lemma}
Under the assumptions of Theorem \ref{thm-main}, there exists $K_{17}>0$ such that
\begin{equation}\label{eq-uniform-u2}
\|u(\cdot,t)\|_{L^2(\Omega)}\le K_{17},
\qquad\text{for all }t\in(0,T_{\max}).
\end{equation}
\end{lemma}
\begin{proof}
Multiplying the first equation of \eqref{eq-model} by $u$ and
integrating by parts, we obtain
\begin{align}\label{eq-cell-L2-identity}
\frac{1}{2}\frac{d}{dt}\int_{\Omega}u^2
=&\int_{\Omega}u\Delta u
  +\int_{\Omega}u\nabla\cdot(u\chi(v)\nabla v)
  +a\int_{\Omega}u^2-b\int_{\Omega}u^3\notag\\
=&-\int_{\Omega}|\nabla u|^2
  -\int_{\Omega}u\chi(v)\nabla u\cdot\nabla v
  +a\int_{\Omega}u^2-b\int_{\Omega}u^3.
\end{align}
Applying Young's and H\"older's inequalities to the chemotactic term in \eqref{eq-cell-L2-identity}, we find
\begin{align}\label{eq-cell-mixed-product}
\left|\int_{\Omega}u\chi(v)\nabla u\cdot\nabla v\right|
\le&\frac{1}{4}\int_{\Omega}|\nabla u|^2
     +\int_{\Omega}u^2\chi(v)^2|\nabla v|^2\notag\\*
\le&\frac{1}{4}\int_{\Omega}|\nabla u|^2
     +\left(\int_{\Omega}\chi(v)^4|\nabla v|^4\right)^{\frac{1}{2}}
       \|u\|_{L^4(\Omega)}^2.
\end{align}
By the two-dimensional Gagliardo--Nirenberg inequality,
\[
\|u\|_{L^4(\Omega)}^2
\le c_1\|u\|_{L^2(\Omega)}\|\nabla u\|_{L^2(\Omega)}
    +c_1\|u\|_{L^2(\Omega)}^2.
\]
Substituting this estimate into \eqref{eq-cell-mixed-product} and applying Young's inequality
to both products, we obtain
\begin{align*}
\left|\int_{\Omega}u\chi(v)\nabla u\cdot\nabla v\right|
\le&\frac{1}{2}\int_{\Omega}|\nabla u|^2
 +c_2\left(1+\int_{\Omega}\chi(v)^4|\nabla v|^4\right)
       \int_{\Omega}u^2.
\end{align*}
Inserting this bound into \eqref{eq-cell-L2-identity} and dropping the nonpositive terms
$-\frac{1}{2}\int_{\Omega}|\nabla u|^2$ and $-b\int_{\Omega}u^3$, we obtain
\begin{equation}\label{eq-logcell}
\frac{d}{dt}\int_{\Omega}u^2
\le c_3\left(1+\int_{\Omega}\chi(v)^4|\nabla v|^4\right)
       \int_{\Omega}u^2.
\end{equation}
By \eqref{eq-u2-window} and \eqref{eq-fourth-power},
\[
\int_t^{t+\tau}\int_{\Omega}u^2\le K_2,\qquad
\int_t^{t+\tau}\left(1+\int_{\Omega}\chi(v)^4|\nabla v|^4\right)ds\le\tau+K_{16}.
\]
Applying Lemma \ref{lem-uniform-gronwall} to \eqref{eq-logcell} with $f=0$,
we obtain $\int_{\Omega}u^2\le c_4$ for $\tau<t<T_{\max}$,
where $c_4=\frac{K_2}{\tau}e^{c_3(\tau+K_{16})}$.
For $0<t\le\tau$, \eqref{eq-uniform-u2} follows from the continuity of $u$ in Lemma \ref{lem-local}.
This completes the proof.
\end{proof}

Having established \eqref{eq-uniform-u2}, we now prove Theorem \ref{thm-main}.

\begin{proof}[Proof of Theorem \ref{thm-main}]
For $0<t<T_{\max}$, the variation-of-constants formula for the second equation of
\eqref{eq-model} gives
\begin{equation}\label{eq-signal-formula}
v(\cdot,t)=e^{t(\Delta-1)}v_0
+\int_0^te^{(t-s)(\Delta-1)}u(\cdot,s)ds.
\end{equation}
Applying the heat semigroup estimates \eqref{eq-heat-Lp}--\eqref{eq-heat-initial-gradient}
and the bound \eqref{eq-uniform-u2} for $u$, we find
\begin{align}
\|v(\cdot,t)\|_{L^\infty(\Omega)}
\le&\|v_0\|_{L^\infty(\Omega)}
+\gamma_1\int_0^t\bigl(1+(t-s)^{-\frac{1}{2}}\bigr)e^{-(t-s)}
\|u(\cdot,s)\|_{L^2(\Omega)}ds
\le c_1,\label{eq-signal-infty}\\
\|\nabla v(\cdot,t)\|_{L^8(\Omega)}
\le&\gamma_3\|\nabla v_0\|_{L^8(\Omega)}
+\gamma_2\int_0^t\bigl(1+(t-s)^{-\frac{7}{8}}\bigr)e^{-(t-s)}
\|u(\cdot,s)\|_{L^2(\Omega)}ds
\le c_2.\notag
\end{align}
With $\|\chi(v)\|_{L^\infty(\Omega)}\le c_3$ by \eqref{eq-signal-infty}, we estimate
the chemotactic term by H\"older's inequality and interpolation to obtain
\begin{align}\label{eq-flux-L4}
\|u\chi(v)\nabla v\|_{L^4(\Omega)}
\le& c_2c_3\|u\|_{L^8(\Omega)}
\le c_2c_3\|u\|_{L^2(\Omega)}^{\frac{1}{4}}
             \|u\|_{L^\infty(\Omega)}^{\frac{3}{4}}\notag\\
\le& c_4\|u\|_{L^\infty(\Omega)}^{\frac{3}{4}}.
\end{align}

For $0<t\le T<T_{\max}$, we apply \eqref{eq-heat-Lp} and \eqref{eq-heat-divergence}
to the variation-of-constants formula for the first equation of \eqref{eq-model}.
Since $e^{t(\Delta-1)}$ preserves positivity and $(a+1)u-bu^2\le(a+1)u$, we obtain
\begin{align*}
\|u(\cdot,t)\|_{L^\infty(\Omega)}
\le&\|u_0\|_{L^\infty(\Omega)}
+\gamma_4\int_0^t\bigl(1+(t-s)^{-\frac{3}{4}}\bigr)e^{-(t-s)}
\|(u\chi(v)\nabla v)(\cdot,s)\|_{L^4(\Omega)}ds\\
&+(a+1)\gamma_1\int_0^t\bigl(1+(t-s)^{-\frac{1}{2}}\bigr)e^{-(t-s)}
\|u(\cdot,s)\|_{L^2(\Omega)}ds.
\end{align*}
Taking the supremum over $[0,T]$ and using \eqref{eq-uniform-u2},
\eqref{eq-flux-L4} and Young's inequality, we obtain
\begin{align*}
\sup_{0\le t\le T}\|u(\cdot,t)\|_{L^\infty(\Omega)}
\le& c_5\left[1+\left(\sup_{0\le t\le T}\|u(\cdot,t)\|_{L^\infty(\Omega)}\right)^{\frac{3}{4}}\right]\\
\le&\frac{1}{2}\sup_{0\le t\le T}\|u(\cdot,t)\|_{L^\infty(\Omega)}+c_6,
\end{align*}
where $c_5,c_6$ are independent of $T$. Letting $T\nearrow T_{\max}$, we obtain
\begin{equation}\label{eq-uniform-u-infty}
\sup_{0\le t<T_{\max}}\|u(\cdot,t)\|_{L^\infty(\Omega)}\le2c_6.
\end{equation}

Finally, we apply \eqref{eq-heat-gradient}--\eqref{eq-heat-initial-gradient}
to \eqref{eq-signal-formula} and use \eqref{eq-uniform-u-infty} to obtain
\[
\|\nabla v(\cdot,t)\|_{L^\infty(\Omega)}
\le\gamma_3\|\nabla v_0\|_{L^\infty(\Omega)}
+\gamma_2\int_0^t\bigl(1+(t-s)^{-\frac{1}{2}}\bigr)e^{-(t-s)}
\|u(\cdot,s)\|_{L^\infty(\Omega)}ds
\le c_7.
\]
The extensibility criterion \eqref{eq-alternative} now gives $T_{\max}=\infty$,
and \eqref{eq-bounded} follows.
This completes the proof of Theorem \ref{thm-main}.
\end{proof}

\end{document}